\documentclass{amsart}
\usepackage{amssymb}
\usepackage{graphics}
\usepackage{amscd}
\usepackage{enumerate}
\newtheorem{thm}{Theorem}[section]
\newtheorem{thm*}{Theorem}[]
\newtheorem{cor}[thm]{Corollary}
\newtheorem{lem}[thm]{Lemma}
\newtheorem{prop}[thm]{Proposition}

\newtheorem{claim}[thm]{Claim}
\theoremstyle{definition}
\newtheorem{defn}[thm]{Definition}

\theoremstyle{remark}
\newtheorem{rem}[thm]{Remark}
\newtheorem{ex}[thm]{Example}
\numberwithin{equation}{section}
\newtheorem{note}[thm]{Notation}

\newcommand{\eps}{\varepsilon}

\begin{document}

\title{The sutured Floer homology polytope}%
\author{Andr\'as Juh\'asz}%
\address{Department of Pure Mathematics and Mathematical Statistics, University of Cambridge, Wilberforce Road, Cambridge, CB3 0WB, UK}%
\email{aij22@dpmms.cam.ac.uk}%

\thanks{Research partially supported by OTKA grant no. T49449.}

\date{\today}%
\begin{abstract}
In this paper, we extend the theory of sutured Floer homology
developed by the author. We first prove an adjunction inequality, and then define a polytope $P(M,\gamma)$ in $H^2(M,\partial M;\mathbb{R})$ that is spanned by the $\text{Spin}^c$-structures which support non-zero Floer homology groups. If $(M,\gamma) \rightsquigarrow (M',\gamma')$ is a taut surface decomposition, then an affine map projects $P(M',\gamma')$ onto a face of $P(M,\gamma);$ moreover, if $H_2(M) = 0,$ then every face
of $P(M,\gamma)$ can be obtained in this way for some surface decomposition. We show that if $(M,\gamma)$ is reduced, horizontally prime, and $H_2(M) = 0,$ then $P(M,\gamma)$ is maximal dimensional in $H^2(M,\partial M;\mathbb{R}).$ This implies that if $\text{rk}(SFH(M,\gamma)) < 2^{k+1},$ then $(M,\gamma)$ has depth at most $2k.$ Moreover, $SFH$ acts as a complexity for balanced sutured manifolds. In particular, the rank of the top term of knot Floer homology bounds the topological
complexity of the knot complement, in addition to simply detecting fibred knots.
\end{abstract}

\maketitle
\section{Introduction}

Heegaard Floer homology is an invariant of closed oriented three-manifolds defined by Ozsv\'ath and Szab\'o in \cite{OSz}. It comes in
fours different flavors: $\widehat{HF},$ $HF^+,$ $HF^-,$ and $HF^{\infty}.$ This was extended to an invariant $HFK$ of knots by Ozsv\'ath and Szab\'o in \cite{OSz3}, and independently by Rasmussen in \cite{Ras}. Later, Ozsv\'ath and Szab\'o \cite{OSz2} generalized $\widehat{HFK}$ to an invariant $\widehat{HFL}$ of links in $S^3.$ For a knot $K$ in $S^3,$ the group $\widehat{HFK}(K)$ splits as a direct sum $\bigoplus_{i, j \in \mathbb{Z}}\widehat{HFK}_j(K,i),$ and has a homological $\mathbb{Z}$-grading. For each $i \in \mathbb{Z},$ the Euler characteristic of $\widehat{HFK}_*(K,i)$ is equal to the $i$-th coefficient $a_i$ of the Alexander-Conway polynomial $\Delta_K(t)$ of $K.$

It is a classical result that for a knot $K$ in $S^3$ the polynomial $\Delta_K(t)$ gives a lower bound on the genus of
$K$ in the following sense: $$g(K) \ge \max \{\, i \in \mathbb{Z} \,\colon\, a_i \neq 0 \,\}.$$ Ozsv\'ath and
Szab\'o in \cite{OSz6} showed that knot Floer homology actually detects the genus of $K:$
$$g(K) = \max \{\, i \in \mathbb{Z} \,\colon\, \widehat{HFK}(K,i) \neq 0 \,\}.$$ The proof of this striking result
uses Gabai's theory of sutured manifolds \cite{Gabai}, \cite{Gabai6}, \cite{Gabai4}, the Eliashberg-Thurston theory
of confoliations \cite{ET}, the contact invariant and cobordism maps in Heegaard Floer homology,
symplectic semi-fillings, and Lefshetz pencils.

The theory of sutured manifolds was developed by Gabai in \cite{Gabai} in order to study the existence of taut foliations on 3-manifolds. Sutured manifolds are oriented 3-manifolds with boundary, together with a set of oriented simple closed curves, called sutures, that divide the boundary into a plus and a minus part. They can be thought of as cobordisms between compact oriented surfaces with boundary. Gabai also defined an operation on sutured manifolds, called sutured manifold decomposition. It consists of cutting the manifold along a properly embedded oriented surface $R,$ and adding one side of $R$ to the plus, and the other side to the minus part of the boundary. He showed that a sutured manifold carries a taut foliation if and only if there is a sequence of decompositions that results in a product sutured manifold (essentially a trivial cobordism). The theory of sutured manifold decompositions was generalized in \cite{convex} to study tight contact structures on 3-manifolds, and was called convex decomposition theory.

In \cite{sutured}, I introduced sutured Floer homology, in short $SFH,$ which is an invariant of balanced sutured manifolds.
$SFH$ is an invariant of three-manifolds with boundary, and generalizes $\widehat{HF}, \widehat{HFK},$ and $\widehat{HFL}.$
The balanced condition is not very restrictive, since in Proposition \ref{prop:4} we show that every open \emph{taut} sutured manifold that has at least one suture on each boundary component is balanced.

$SFH$ was used in \cite{decomposition} to give a more elegant and direct proof of the fact that knot Floer homology detects
the genus of a knot. That proof only relies on Gabai's theory of sutured manifolds and the following two results.
First, if $R$ is a Seifert surface of a knot $K$ in $S^3,$ then $$\widehat{HFK}(K,g(R)) \cong SFH(S^3(R)),$$
where $S^3(R)$ is the sutured manifold complementary to $R,$ see \cite[Theorem 1.5]{decomposition}. Secondly, by \cite[Theorem 1.3]{decomposition}, if we decompose a sutured manifold $(M,\gamma)$
along a ``nice" surface and get the sutured manifold $(M',\gamma'),$ then $SFH(M',\gamma')$ is a direct summand of
$SFH(M,\gamma).$ We will refer to this as ``the decomposition formula". If $R$ is of minimal genus, then $S^3(R)$ is taut, so by \cite{Gabai} there is a sequence of nice decompositions
that ends in a product. The $SFH$ of a product is $\mathbb{Z},$ so the decomposition formula implies that $SFH(S^3(R))$ contains a $\mathbb{Z}$ direct summand.

For a Seifert surface $R,$ even though $SFH(S^3(R))$ is isomorphic to the top term of knot Floer homology, it carries
an extra $\text{Spin}^c$-grading. Note that for a sutured manifold $(M,\gamma),$ the set of $\text{Spin}^c$-structures
$\text{Spin}^c(M,\gamma)$ is an affine space over $H^2(M,\partial M) \cong H_1(M).$ In the present paper, we study this
extra grading on $SFH,$ and how it behaves under sutured manifold decompositions. Using our results, we show that the
top term of knot Floer homology carries deep topological information about the knot complement. In particular, we have the following,
which is a special case of Corollary \ref{cor:7}.

\begin{thm*} \label{thm:depth}
Suppose that $K$ is a knot in $S^3,$ and
$$\text{rk}\left(\widehat{HFK}(K,g(K))\right) < 2^{k+1}.$$ Then the sutured manifold $S^3(K)$ complementary
to $K$ has depth $d(Y(K)) \le 2k+1.$ In particular, if $k=0,$ then $K$ is fibred.
\end{thm*}

Here the depth of a sutured manifold is the minimal number of decompositions needed to get a product sutured manifold.

Ozsv\'ath and Szab\'o conjectured that knot Floer homology detects fibred knots in the sense that
$\widehat{HFK}(K,g(K)) \cong \mathbb{Z}$ if and only if $K$ is fibred. This was proved by Ghiggini in \cite{Ghiggini} for genus one knots, and proceeds along the lines of the Ozsv\'ath-Szab\'o proof of
$\widehat{HFK}(K,g(K)) \neq 0,$ using deep symplectic and contact topology. Building on Ghiggini's work, Ni proposed a proof of the general case in \cite{fibred}, using an alternative version of sutured Floer homology (without the $\text{Spin}^c$-grading), and a restricted version of the decomposition formula for horizontal surfaces and separating product annuli. Shortly after this, in \cite{decomposition} I presented a more direct proof of the fibred knot conjecture, only using $SFH$ and the general decomposition formula. This starts out with an observation of Gabai \cite{Gabai2} that a knot $K$ is fibred if and only if $S^3(R)$ is a product sutured manifold, where $R$ is a minimal genus Seifert surface for $K.$ So the problem can be reduced to the question whether $SFH$ detects product sutured manifolds.  Later, it turned out that the last part of the proof in \cite{fibred} had a gap due to an incorrect reference to \cite{cpr} concerning characteristic product regions. In \cite{decomposition}, I borrowed Ni's last argument to conclude my proof, so \cite{decomposition} has the same gap.
Ni filled in this gap in \cite{corrigendum}. In the present paper, I correct and generalize \cite{decomposition} by eliminating
the use of characteristic product regions. My approach is completely different from that of \cite{corrigendum}. Instead, I only use reduced sutured manifolds, ones in which every product annulus is parallel to a suture. Since in \cite{decomposition} I also proved the decomposition formula for non-separating product annuli, it is enough to work with reduced sutured manifolds. The introduction of the
$SFH$ polytope makes the proof very transparent, and makes it possible to get a much sharper result, namely Theorem \ref{thm:depth}.

For a sutured manifold $(M,\gamma),$ the $\text{Spin}^c$-structures that support non-zero Floer homology groups span a polytope $P(M,\gamma)$ in $H^2(M,\partial M; \mathbb{R}).$ This polytope is well defined up to translations. A major tool in this paper is the adjunction inequality, Theorem \ref{thm:1}.

\begin{thm*}[Adjunction Inequality]
Suppose that the sutured manifold $(M,\gamma)$ is strongly balanced,
and fix a trivialization $t \in T(M,\gamma).$ Let $S \subset M$ be a
nice decomposing surface. If a $\text{Spin}^c$-structure
$\mathfrak{s} \in \text{Spin}^c(M,\gamma)$ satisfies $$\langle\,
c_1(\mathfrak{s},t),[S]\,\rangle < c(S,t),$$ then $SFH(M,\gamma,\mathfrak{s}) =
0.$
\end{thm*}

Here $c(S,t)$ is a purely topological quantity, and we show that the above inequality can be rearranged to get a Thurston-Bennequin
type inequality. In Theorem \ref{thm:3}, we use the adjunction inequality to extend the decomposition formula \cite[Theorem 3.11]{decomposition} to disconnected decomposing surfaces.

In Proposition \ref{prop:1} and Corollary \ref{cor:4}, we establish a relationship between decompositions of $(M,\gamma)$ and faces of $P(M,\gamma),$ and show that if $H_2(M)=0,$ then every face of $P(M,\gamma)$ corresponds to a well-groomed surface decomposition. More concretely, Theorem \ref{thm:5} implies that if $(M,\gamma) \rightsquigarrow^S (M',\gamma')$ is a taut surface decomposition, then there is an affine map from $H^2(M',\partial M'; \mathbb{R})$ to $H^2(M,\partial M;\mathbb{R})$ which projects $P(M',\gamma')$ onto a face of $P(M,\gamma).$ This map is a translate of the dual of the map $H_1(M') \to H_1(M)$ induced by the embedding $M' \hookrightarrow M.$ If $S$ is a disk, then this projection is actually an isomorphism. So we see how $\text{Spin}^c$-structures split under surface decompositions. This, for example, implies a result of Gabai \cite{Gabai7} that if a sutured manifold is disk decomposable, then it can be decomposed into a product using a single (not necessarily connected) surface, and if $\gamma$ is connected, then it carries a taut foliation of depth at most one.

From now on, we are going to suppose that $(M,\gamma)$ is a taut balanced sutured manifold which satisfies the condition $H_2(M) = 0.$ The condition $H_2(M) = 0$ is not very restrictive, since it is satisfied by any sutured manifold complementary to a connected surface in a rational homology 3-sphere; furthermore, it is preserved by nice surface decompositions. And the most studied sutured manifolds are exactly the ones which are complementary to a Seifert surface of a knot or a link.

Theorem \ref{thm:2} is one of the main results of this paper.

\begin{thm*} \label{thm:dim}
Suppose that $H_2(M) = 0,$ and the sutured manifold $(M,\gamma)$ is
balanced, taut, reduced and horizontally prime. Then $$\dim P(M,\gamma) = \dim H^2(M,\partial M;\mathbb{R}) = b_1(M) = b_1(\partial M)/2.$$
In particular,
$$\text{rk}(SFH(M,\gamma)) \ge b_1(\partial M)/2+1.$$
\end{thm*}

This result fills in the gap in \cite{decomposition}, and makes further generalizations possible. Using this, we prove Proposition \ref{prop:6}:

\begin{thm*} \label{thm:depth2}
Suppose that $(M,\gamma)$ is a taut balanced sutured manifold such
that $H_2(M) = 0$ and $\text{rk}(SFH(M,\gamma)) < 2^{k+1}$ for some
integer $k \ge 0.$ Then the depth of $(M,\gamma)$ is at most $2k.$
\end{thm*}

The proof proceeds by induction on $k.$ One has to first decompose $(M,\gamma)$ along a maximal set of product annuli to make it reduced. Then $P(M,\gamma)$ becomes maximal dimensional in its ambient space. There is a $\text{Spin}^c$-structure $\mathfrak{s}$ that is a vertex of $P(M,\gamma),$ and such that $\text{rk}(SFH(M,\gamma,\mathfrak{s})) < 2^k.$ Furthermore, we saw that there is a decomposition $(M,\gamma) \rightsquigarrow^S (M',\gamma')$ such that $SFH(M',\gamma') \cong SFH(M,\gamma,\mathfrak{s}).$ So we can apply the induction hypotheses to $(M',\gamma')$ to see that it has depth at most $2k-2.$ In particular, this illustrates how
the rank of $SFH$ can be used to measure the complexity of balanced sutured manifolds, and to perform inductive proofs using it.
It is worth comparing it to the complexity defined by Gabai in \cite{Gabai} to show the existence of sutured manifold hierarchies.

Theorem \ref{thm:depth2} implies Theorem \ref{thm:depth},
since if we decompose the knot complement $S^3(K)$ along a minimal genus Seifert surface $R,$ then we get the taut balanced
sutured manifold $S^3(R)$ with $\text{rk}(SFH(M,\gamma)) < 2^{k+1}.$

The sutured manifold $(M,\gamma)$ constructed in Example \ref{ex:1} has the following surprising property: the polytope $P(M,\gamma)$ consists of a single point, even though $(M,\gamma)$ is horizontally prime and is not a product. This example also illustrates that decompositions along product annuli can change the sutured Floer homology polytope, while the rank of $SFH$ remains unchanged.
So we always need to make our sutured manifold reduced by cutting along product annuli before we can decrease the rank of $SFH$
by another sutured manifold decomposition.

In Section \ref{sect:7}, we define a function $y$ on $H_2(M,\partial M;\mathbb{R})$ that is a semi-norm, except that $y(c)$ and $y(-c)$ might be different. The dual unit norm polytope of $y$ is exactly $-P(M,\gamma).$ Moreover, $y$ is non-degenerate if $(M,\gamma)$ is reduced and horizontally prime. Using this norm, we can view Corollary \ref{cor:1} as an extension of a theorem of Ozsv\'ath and Szab\'o \cite{OSz7} that link Floer homology detects the Thurston norm of the link complement. We will prove in \cite{decateg} that
if we symmetrize $y,$ we get a semi-norm that gives a lower bound on the semi-norm defined by Scharlemann in \cite{Scharlemann},
but is different from it.

Finally, we compute the sutured Floer homology of any sutured manifold $(M,\gamma)$ such that $M \approx S^1 \times D^2.$ This illustrates some of the techniques developed in this paper, and will be used in future computations. For further examples of $P(M,\gamma),$ including whole families where $M$ is a genus two handlebody, we refer the reader to \cite{decateg}. In \cite{HJS}, we show how the sutured Floer homology polytope can be used to distinguish Seifert surfaces up to isotopy.

\section*{Acknowledgement}

I am extremely grateful for the guidance of David Gabai during the course of this work.
His insight on sutured manifold theory and the class he gave in the academic year 2007/08 at Princeton University on the theory of foliations were invaluable. I would also like to thank Zolt\'an Szab\'o and Paolo Ghiggini for several helpful discussions. The paper was rewritten in July 2008 at the Institut des Hautes \'Etudes Scientifiques, and in September 2008 at the R\'enyi Institute, where the author was supported by the BudAlgGeo (Algebraic Geometry) project, in the framework of the European Community's ``Structuring the European Research Area" programme.

\section{Sutured manifolds}
To get an in depth introduction to the theory of sutured manifolds and surface decompositions, we recommend reading Gabai's original papers \cite{Gabai, Gabai6, Gabai4}. For the reader's convenience, let us review the most important definitions and results here.

\begin{note}
Throughout this paper, we are going to use the following notations. If
$K$ is a submanifold of the manifold $M,$ then $N(K)$ denotes a
regular neighborhood of $K$ in $M$ and $[K]$ is the homology class represented by $K.$
If $A$ is a set, then $|A|$ is the cardinality of $A.$
If $X$ is a topological space, then $|X|$ is the number of components of $X.$
\end{note}

\begin{defn}
A \emph{sutured manifold} $(M,\gamma)$ is a compact oriented
3-manifold $M$ with boundary together with a set $\gamma \subset
\partial M$ of pairwise disjoint annuli $A(\gamma)$ and tori
$T(\gamma).$ Furthermore, the interior of each component of
$A(\gamma)$ contains a \emph{suture}, i.e., a homologically
nontrivial oriented simple closed curve. We denote the union of the
sutures by $s(\gamma).$

Finally every component of $R(\gamma)=\partial M \setminus
\text{Int}(\gamma)$ is oriented. Define $R_+(\gamma)$ (or
$R_-(\gamma)$) to be those components of $\partial M \setminus
\text{Int}(\gamma)$ whose normal vectors point out of (into) $M$.
The orientation on $R(\gamma)$ must be coherent with respect to
$s(\gamma),$ i.e., if $\delta$ is a component of $\partial
R(\gamma)$ and is given the boundary orientation, then $\delta$ must
represent the same homology class in $H_1(\gamma)$ as some suture.
\end{defn}

\begin{defn}
A sutured manifold $(M,\gamma)$ is called \emph{balanced} if M has
no closed components, $\chi(R_+(\gamma))=\chi(R_-(\gamma)),$ and the
map $\pi_0(A(\gamma)) \to \pi_0(\partial M)$ is surjective.
\end{defn}

\begin{defn}
A sutured manifold $(M,\gamma)$ is \emph{taut} if $M$ is irreducible
and $R(\gamma)$ is incompressible and Thurston norm minimizing in $H_2(M,\gamma).$
\end{defn}

\begin{defn} \label{defn:8}
Let $(M, \gamma)$ be a sutured manifold. A \emph{decomposing
surface} is a properly embedded oriented surface $S$ in $M$ such
that no component of $\partial S$ bounds a disk in $R(\gamma)$ and no component of
$S$ is a disk $D$ with $\partial D \subset R(\gamma).$ Moreover,
for every component $\lambda$ of $S \cap \gamma$ one of (1)-(3) holds:
\begin{enumerate}
\item $\lambda$ is a properly embedded non-separating arc in $\gamma$
such that $|\lambda \cap s(\gamma)| = 1.$
\item $\lambda$ is a simple closed curve in an annular component $A$
of $\gamma$ in the same homology class as $A \cap s(\gamma).$
\item $\lambda$ is a homotopically nontrivial curve in a torus
component $T$ of $\gamma,$ and if $\delta$ is another component of
$T \cap S,$ then $\lambda$ and $\delta$ represent the same homology
class in $H_1(T).$
\end{enumerate}

Then $S$ defines a \emph{sutured manifold decomposition}
$$(M, \gamma)\rightsquigarrow^{S} (M', \gamma'),$$ where $M' = M
\setminus \text{Int}(N(S))$ and $$\gamma' = (\gamma \cap M') \cup
N(S'_+ \cap R_-(\gamma)) \cup N(S'_- \cap R_+(\gamma)), $$
$$R_+(\gamma') = ((R_+(\gamma) \cap M') \cup S'_+) \setminus
\text{Int}(\gamma'),$$
$$R_-(\gamma') = ((R_-(\gamma) \cap M') \cup S'_-) \setminus
\text{Int}(\gamma'),$$ where $S'_+$ ($S'_-$) is the component of
$\partial N(S) \cap M'$ whose normal vector points out of (into)
$M'.$
\end{defn}

\begin{defn}
If $(M,\gamma)$ is a balanced sutured manifold then a surface decomposition $(M,\gamma) \rightsquigarrow^S (M',\gamma')$ is called \emph{groomed} if for each component $V$ of $R(\gamma)$ either $S \cap V$ is a union of parallel, coherently oriented, nonseparating
closed curves or $S \cap V$ is a union of arcs such that for each component $\delta$ of
$\partial V$ we have $|\delta \cap \partial S| = |\langle\, \delta \cap \partial S \,\rangle|.$

A surface decomposition is called \emph{well groomed} if for each component $V$ of
$R(\gamma)$ it holds that $S \cap V$ is a union of parallel, coherently oriented, nonseparating closed curves or arcs.
\end{defn}

The following definition is motivated by \cite[Lemma 3.8]{Gabai}.

\begin{defn}
Let $(M,\gamma)$ be a balanced sutured manifold. We say that a class $z \in H_2(M,\partial M)$ is \emph{well groomed}
if $\partial z \neq 0$ in $H_1(\partial M)$ and the following hold.

(1) For each non-planar component $V$ of $R(\gamma)$ and each component $\lambda$ of $\partial V$ we have $\langle\,
z,\lambda\,\rangle = 0.$

(2) For each planar component $V$ of $R(\gamma)$ there exist at most two components $\lambda_1$ and $\lambda_2$ of $\partial V$
such that $\langle\, z,\lambda_i\,\rangle \neq 0$ for $i = 1,2.$
\end{defn}

Note that $z \in H_2(M,\partial M)$ is well groomed if and only if $-z$ is well groomed. Using this terminology \cite[Lemma 3.8]{Gabai} can be stated as follows.

\begin{lem} \label{lem:5}
Let $(M,\gamma)$ be a balanced sutured manifold. Then there exists a well groomed class in $H_2(M,\partial M).$
\end{lem}

\begin{lem} \label{lem:6}
Let $(M,\gamma)$ be a taut balanced sutured manifold and $z \in H_2(M,\partial M)$ a well groomed homology class.
Then there is a well groomed surface decomposition $(M,\gamma) \rightsquigarrow^S (M',\gamma')$ such that $[S] = z$
and $(M',\gamma')$ is taut.
\end{lem}

\begin{proof}
This lemma follows from the last argument in the proof of \cite[Lemma 3.4]{Gabai4} which goes as follows. By \cite[Lemma 0.7]{Gabai6}
there is a groomed surface $S$ which gives a taut decomposition and $[S] = z.$ Now \cite[Lemma 3.9]{Gabai} implies that
if $V$ is a component of $R(\gamma)$ then $S \cap V$ is homologous to a set of parallel curves. Finally, an application of \cite[Lemma 0.6]{Gabai6} yields the desired well groomed surface.
\end{proof}

\begin{defn}
Let $(M,\gamma)$ be a sutured manifold. A \emph{product annulus} in $(M,\gamma)$ is an annulus $A$ properly embedded in $M$
such that $\partial A \subset R(\gamma),$ $\partial A \cap R_+(\gamma) \neq \emptyset,$
and $\partial A \cap R_-(\gamma) \neq \emptyset.$ A \emph{product disk} is a disk $D$
properly embedded in $M$ such that $\partial D \cap \gamma$
consists of two essential arcs in $\gamma.$ Product disks and product annuli detect
where a sutured manifold is locally a product. $(M,\gamma)$ is a \emph{product sutured manifold}
if $M = R \times I,$ $\gamma = \partial R \times I,$ $R_+(\gamma) = R \times \{1\},$ and
$R_-(\gamma) = R \times \{0\}.$
\end{defn}

\begin{lem} \label{lem:7}
Let $(M,\gamma)$ be a sutured manifold such that $M$ is irreducible and $R(\gamma)$ is incompressible. If $A$ is a compressible product annulus in $(M,\gamma),$ then there is a
cylinder $D^2 \times I \subset M$ such that $A = \partial D^2 \times I$ and $D^2 \times \partial I \subset R(\gamma).$
\end{lem}

\begin{proof}
Suppose that a simple closed curve in $A$ bounds a disk in $M$ but it does not bound a disk in $A.$ Then both $a_- = A \cap R_-(\gamma)$ and $a_+ = A \cap R_+(\gamma)$ bound disks in $M.$ Since $R(\gamma)$ is incompressible, both $a_-$ and $a_+$ bound disks $D_-$ and $D_+$ in $R(\gamma),$ respectively.
But $M$ is irreducible, hence the embedded sphere $D_- \cup A \cup D_+$ bounds a 3-ball in $M.$
This 3-ball can be identified with $D^2 \times I$ such that it satisfies the stated properties.
\end{proof}

\begin{defn}
We say that a balanced sutured manifold $(M,\gamma)$ is \emph{reduced} if every incompressible product annulus $A$ in $(M,\gamma)$ is ambient isotopic to a component of $\gamma$ such that $\partial A$ stays in $R(\gamma)$ throughout. Call a product disk $D$ \emph{inessential} if there is an ambient isotopy of $D$ into $\gamma$ which fixes $D \cap \gamma,$ and \emph{essential} otherwise.
\end{defn}

\begin{lem} \label{lem:8}
Let $(M,\gamma)$ be a \emph{reduced} sutured manifold such that $M$ is irreducible and $R(\gamma)$ is incompressible. Then exactly one of the following holds.
\begin{enumerate}
\item Every product disk in $(M,\gamma)$ is \emph{inessential}.
\item $(M,\gamma)$ is homeomorphic to $(\Sigma \times I, \partial \Sigma \times I),$ where $\Sigma$ is a sphere with either two or three open disks removed.
\end{enumerate}
\end{lem}

\begin{proof}
In this proof we implicitly use the following observation several times. Suppose that the product annulus $A \subset M \setminus \gamma$ is ambient isotopic to a component $\gamma'$ of $\gamma$ such that $\partial A$ stays in $R(\gamma)$ throughout. Then there is a submanifold $C \times I$ inside $M$ such that $C$ is an annulus, $\partial C \times I = A \cup \gamma'$ and $C \times \partial I \subset R(\gamma).$

Suppose that $(M,\gamma)$ contains an essential product disk $D.$ We distinguish two cases depending on whether the two arcs of $D \cap \gamma$ lie in the same component of $\gamma.$
First suppose that there is a single component $\gamma_0$ of $\gamma$ which contains $D \cap \gamma.$ Denote the components of $\partial N(\gamma_0 \cup D) \setminus \gamma_0$ by $A_1$ and $A_2.$ Then the product annuli $A_1$ and $A_2$ both have to be incompressible, otherwise by Lemma \ref{lem:7} the product disk $D$ would by inessential. Furthermore, neither $A_1$ nor $A_2$ can be ambient isotopic to $\gamma_0,$ else again $D$ would be inessential. Thus there are components $\gamma_1$ and $\gamma_2$ of $\gamma,$ both distinct from $\gamma_0$ and from each other, such that $A_i$ is ambient isotopic to $\gamma_i$ for $i = 1,2.$ It follows that $(M,\gamma)$ is the product $(\Sigma \times I, \partial \Sigma \times I),$ where $\Sigma$ is a sphere with three open disks removed.

Now suppose that there are components $\gamma_0$ and $\gamma_1$ of $\gamma$ such that $D \cap \gamma_i \neq \emptyset$ for $i=0,1.$ Let $A = \partial N(\gamma_0 \cup \gamma_1 \cup D) \setminus (\gamma_0 \cup \gamma_1).$ If the product annulus $A$ is compressible then by Lemma \ref{lem:7} we are in case (2) with $\Sigma$ being a sphere with two open disks removed. Otherwise $A$ is ambient isotopic to a component $\gamma_2$ of $\gamma$ different from $\gamma_0$ and $\gamma_1,$ and hence we are again in case (2) with $\Sigma$ being a sphere with three open disks removed.

On the other hand, if (2) holds, then $(M,\gamma)$ is reduced, but it contains an essential product disk.
\end{proof}

Next we recall \cite[Proposition V.1.6]{Jaco-Shalen}.

\begin{prop} \label{prop:13}
For each compact, irreducible 3-manifold pair $(M,T),$ there is a number $h(M,T)$ with the following property. Let $W \subset M$ be a two-sided, incompressible surface having more than $h(M,T)$ components and such that $\partial W \subset T.$ Then either
\begin{enumerate}
\item $W$ has a $T$-parallel component, or
\item $W$ has two components which are parallel in $(M,T).$
\end{enumerate}
\end{prop}

\begin{prop} \label{prop:14}
Let $(M,\gamma)$ be a sutured manifold such that $M$ is irreducible and $R(\gamma)$ is incompressible. Then there is a decomposition $(M,\gamma) \rightsquigarrow^A (M',\gamma')$ such that $A$ is a union of pairwise disjoint incompressible product annuli and $(M',\gamma')$ is reduced.
\end{prop}

\begin{proof}
Using the terminology of \cite{Jaco-Shalen} the 3-manifold pair $(M,R(\gamma))$ is irreducible, thus we can apply Proposition \ref{prop:13} to get a number $h(M,R(\gamma)).$ Note that a product annulus cannot be $R(\gamma)$-parallel since its two boundary components lie in different components of $R(\gamma).$ So we can recursively construct a maximal set of pairwise disjoint incompressible product annuli $A_1,\dots,A_n$ such that for $A = A_1 \cup \dots \cup A_n$ no two components of $\gamma \cup A$ are parallel in $(M,R(\gamma)).$ Indeed, $n \le h(M,R(\gamma)) - |\gamma|$ for such an $A,$ thus the recursion has to terminate in finitely many steps.

Decomposing $(M,\gamma)$ along a maximal $A$ the resulting $(M',\gamma')$ is reduced. Indeed, an incompressible product annulus $C'$ in $(M',\gamma')$ which is not parallel to $\gamma'$ gives rise to a product annulus $C$ in $(M,\gamma)$ which is not parallel to any component of $\gamma \cup A.$ We show that $C$ is incompressible. Indeed, if $C$ was compressible, then $\partial C \cap R_+(\gamma)$ would bound a disk $D$ in $R_+(\gamma).$ Since $C'$ is incompressible $\partial A \cap D \neq \emptyset,$ thus $A$ would also be compressible, a contradiction. But the existence of such a $C$ contradicts the maximality of $A.$
\end{proof}

\begin{defn}
Let $(M,\gamma)$ be a balanced sutured manifold. A decomposing
surface $S \subset M$ is called a \emph{horizontal surface} if
\begin{enumerate}[i{)}]
\item $S$ is open and incompressible,
\item $\partial S \subset \gamma$ and $\partial S$ is isotopic to $\partial R_+(\gamma),$
\item $[S] = [R_+(\gamma)]$ in $H_2(M,\gamma),$
\item $\chi(S) = \chi(R_+(\gamma)).$
\end{enumerate}
We say that $(M,\gamma)$ is \emph{horizontally prime} if every
horizontal surface in $(M,\gamma)$ is parallel to either
$R_+(\gamma)$ or $R_-(\gamma).$
\end{defn}

\begin{prop} \label{prop:15}
Let $(M,\gamma)$ be a balanced sutured manifold. Then there is a surface decomposition $(M,\gamma) \rightsquigarrow^H (M',\gamma')$ such that every component of $H$ is a horizontal surface and $(M',\gamma')$ is horizontally prime.
\end{prop}

\begin{proof}
Apply Proposition \ref{prop:13} to the 3-manifold pair $(M,\gamma).$
\end{proof}

\section{Sutured Floer homology and $\text{Spin}^c$ structures}

Sutured Floer homology is an invariant of balanced sutured manifolds defined in \cite{sutured}.
It is constructed in a way analogous to ordinary Heegaard Floer homology.

\begin{defn}
A \emph{sutured Heegaard diagram} is a tuple $( \Sigma,
\boldsymbol{\alpha}, \boldsymbol{\beta}),$ where $\Sigma$ is a
compact oriented surface with boundary and $\boldsymbol{\alpha}$ and
$\boldsymbol{\beta}$ are two sets of pairwise disjoint simple closed
curves in $\text{Int}(\Sigma).$
\end{defn}

Every sutured Heegaard diagram $( \Sigma, \boldsymbol{\alpha},
\boldsymbol{\beta})$ uniquely \emph{defines} a sutured manifold $(M,
\gamma)$ using the following construction. Suppose that
$\boldsymbol{\alpha}=\{\,\alpha_1,\dots,\alpha_m\,\}$ and
$\boldsymbol{\beta}=\{\,\beta_1,\dots,\beta_n\,\}.$ Let $M$ be the
3-manifold obtained from $\Sigma \times I$ by attaching
3-dimensional 2-handles along the curves $\alpha_i \times \{0\}$ and
$\beta_j \times \{1\}$ for $i=1, \dots, m$ and $j=1, \dots, n.$ The
sutures are defined by taking $\gamma =
\partial \Sigma \times I$ and $s(\gamma)=
\partial \Sigma \times \{1/2\}.$

Let $(\Sigma, \boldsymbol{\alpha},\boldsymbol{\beta})$ be an admissible sutured Heegaard diagram defining a balanced sutured manifold $(M,\gamma).$ (For the definition of admissibility see \cite[Definition 3.11]{sutured}, it means that every non-zero periodic domain has both positive and negative coefficients.) Then $|\boldsymbol{\alpha}| = |\boldsymbol{\beta}|,$ denote this number by  $d.$ After an appropriate choice of a generic almost complex and a symplectic structure on $\text{Sym}^d(\Sigma),$
we can apply the Lagrangian Floer homology machinery to the Lagrangian submanifolds $\mathbb{T}_{\alpha} = \alpha_1 \times \dots \times \alpha_d$ and $\mathbb{T}_{\beta} = \beta_1 \times \dots \times \beta_d$ of $\text{Sym}^d(\Sigma).$ This way we obtain a chain complex whose homology $SFH(M,\gamma)$ depends only on the homeomorphism type of $(M,\gamma).$ For the details
see \cite{sutured}. Now we recall \cite[Corollary 3.12]{sutured}.

\begin{lem} \label{lem:3}
Let $(M,\gamma)$ be a balanced sutured manifold such that $H_2(M) = 0.$ Then
every balanced diagram defining $(M,\gamma)$ is admissible.
\end{lem}

Next we review the definition of a $\text{Spin}^c$ structure on a
balanced sutured manifold $(M,\gamma),$ which was introduced in
\cite{sutured}, also see \cite{decomposition}.
Note that in a balanced sutured manifold none of the
sutures are tori. Fix a Riemannian metric on $M.$

\begin{note}
Let $v_0$ be a nowhere vanishing vector field along $\partial M$
that points into $M$ along $R_-(\gamma),$ points out of $M$ along
$R_+(\gamma),$ and on $\gamma$ it is the gradient of a height
function $s(\gamma) \times I \to I.$ The space of such vector fields
is contractible.
\end{note}

\begin{defn}
Let $v$ and $w$ be nowhere vanishing vector fields on $M$ that agree
with $v_0$ on $\partial M.$ We say that $v$ and $w$ are
\emph{homologous} if there is an open ball $B \subset \text{Int}(M)$
such that $v|(M \setminus B)$ is homotopic to $w|(M \setminus B)$
through nowhere vanishing vector fields rel $\partial M.$ We define
$\text{Spin}^c(M, \gamma)$ to be the set of homology classes of
nowhere vanishing vector fields $v$ on $M$ such that $v|\partial M =
v_0.$
\end{defn}

\begin{prop}
Let $(M,\gamma)$ be a sutured manifold such that $M$ is open and the map $\pi_0(\gamma) \to \pi_0(\partial M)$ is surjective. Then $\text{Spin}^c(M,\gamma) \neq \emptyset$ if and only if every component of $M$ is balanced.
\end{prop}

\begin{proof}
It is sufficient to prove the proposition for $M$ connected. First suppose that $\text{Spin}^c(M,\gamma) \neq \emptyset,$ and let $\mathfrak{s}$ be an arbitrary element. Recall from \cite[Definition 4.4]{sutured} that $c_1(\mathfrak{s})$ is defined as the Euler class of the oriented 2-plane field $v^{\perp},$ where $v$ is an arbitrary vector field representing $\mathfrak{s}.$ Let $i \colon \partial M \hookrightarrow M$ denote the embedding and let $\delta = c_1(v_0^{\perp}).$ Then $i^*(c_1(\mathfrak{s})) = \delta,$ thus $$\langle \delta, [\partial M] \,\rangle = \langle\, i^*(c_1(\mathfrak{s})), [\partial M] \,\rangle = \langle c_1(\mathfrak{s}), i_*([\partial M]) \rangle = 0$$ since the cycle $\partial M$ represents zero in $H_2(M;\mathbb{Z}).$ On the other hand, $v_0^{\perp}|R_+(\gamma) = TR_+(\gamma)$ and $v_0^{\perp}|R_-(\gamma) = -TR_-(\gamma),$ so
$$\langle\, \delta, [\partial M] \,\rangle = \chi(R_+(\gamma)) - \chi(R_-(\gamma)).$$ Thus $(M,\gamma)$ is balanced.

Now suppose that $(M,\gamma)$ is balanced. Let $f$ be a Morse function as in the proof of
\cite[Proposition 2.13]{sutured}. Then the
vector field $\text{grad}(f)|\partial M = v_0,$ the number $d$ of
index 1 and 2 critical points of $f$ agree, and $f$ has no index 0
or 3 critical points. Choose $d$ pairwise disjoint balls in $M,$
each containing exactly one index 1 and one index 2 critical point
of $f.$ Then we can modify $\text{grad}(f)$ on these balls so that
we obtain a nowhere zero vector field on $M$ such that $v|\partial M =
v_0.$ This shows that $\text{Spin}^c(M,\gamma) \neq \emptyset.$
\end{proof}

\begin{defn}
Let $(M, \gamma)$ be a balanced sutured manifold and $(\Sigma,
\boldsymbol{\alpha},\boldsymbol{\beta})$ a balanced diagram defining
it. To each $\mathbf{x} \in \mathbb{T}_{\alpha} \cap
\mathbb{T}_{\beta}$ we assign a $\text{Spin}^c$ structure
$\mathfrak{s}(\mathbf{x}) \in \text{Spin}^c(M, \gamma)$ as follows.
Choose a Morse function $f$ on $M$ compatible with the given
balanced diagram $(\Sigma, \boldsymbol{\alpha},
\boldsymbol{\beta}).$ Then $\mathbf{x}$ corresponds to a
multi-trajectory $\gamma_{\mathbf{x}}$ of $\text{grad}(f)$
connecting the index one and two critical points of $f$. In a
regular neighborhood $N(\gamma_{\mathbf{x}})$ we can modify
$\text{grad}(f)$ to obtain a nowhere vanishing vector field $v$ on
$M$ such that $v|\partial M = v_0.$ We define
$\mathfrak{s}(\mathbf{x})$ to be the homology class of this vector
field $v.$
\end{defn}

\begin{defn}
We call a sutured manifold $(M,\gamma)$ \emph{strongly balanced} if
for every component $F$ of $\partial M$ the equality $\chi(F \cap
R_+(\gamma)) = \chi(F \cap R_-(\gamma))$ holds.
\end{defn}

The following is \cite[Proposition 3.4]{decomposition}.

\begin{prop}
The vector bundle $v_0^{\perp}$ over $\partial M$ is trivial if and
only if $(M,\gamma)$ is strongly balanced.
\end{prop}

\begin{note}
If the sutured manifold $(M,\gamma)$ is strongly balanced, then let
$T(M,\gamma)$ denote the set of trivializations of $v_0^{\perp}.$
\end{note}

\begin{defn}
Suppose that $(M,\gamma)$ is a strongly balanced sutured manifold.
Let $t \in T(M,\gamma) $ and $\mathfrak{s}
\in \text{Spin}^c(M,\gamma).$ Then we define $$c_1(\mathfrak{s},t)
\in H^2(M, \partial M; \mathbb{Z})$$ to be the relative Euler class
of the vector bundle $v^{\perp}$ with respect to the trivialization
$t.$ In other words, $c_1(\mathfrak{s},t)$ is the obstruction to
extending $t$ from $\partial M$ to a trivialization of $v^{\perp}$
over $M.$
\end{defn}

\begin{lem} \label{lem:2}
Suppose that $(M,\gamma)$ is a strongly balanced sutured manifold.
Let $d \colon H^1(\partial M) \to H^2(M,\partial M)$ be the co-boundary
map in the cohomology long exact sequence of the pair $(M,\partial
M).$ If $\mathfrak{s} \in \text{Spin}^c(M,\gamma)$ and $t_1,t_2 \in
T(M,\gamma),$ then $$c_1(\mathfrak{s},t_1) - c_1(\mathfrak{s},t_2) =
d(t_1 - t_2).$$
\end{lem}

\begin{proof}
Fix a nowhere vanishing vector field $v$ on $M$ representing the $\text{Spin}^c$-structure
$\mathfrak{s}$ and also fix a triangulation of $M.$ The circle bundle $Sv^{\perp}$ over $M$ will be denoted by $E.$ A trivialization $t \in T(M,\gamma)$ can be considered to be a section of $E|\partial M,$ and by definition $c_1(\mathfrak{s},t)$ is the obstruction to extending $t$ from $\partial M$ to $M.$ More precisely, choose an arbitrary extension of $t$ to the one-skeleton of $M.$ Then
the value of a co-cycle $o(E,t)$ representing $c_1(\mathfrak{s},t)$ on a two-simplex $\Delta$ is the homotopy class of $t|\partial \Delta$ in $\pi_1(S^1) \cong \mathbb{Z}$ obtained after trivializing $E|\Delta.$

Given the sections $t_1$ and $t_2$ of $E$ over $\partial M,$ we can homotope them to coincide on the zero-skeleton of $\partial M.$ Then we can choose a common extension of $t_1$ and $t_2$ to $\text{sk}_1(M) \setminus \text{sk}_1(\partial M).$ The cohomology class $t_1-t_2 \in H^1(\partial M;\mathbb{Z})$ is represented by the co-cycle $o(t_1,t_2).$ The value of $o(t_1,t_2)$ on an edge $\epsilon$ of $\text{sk}_1(\partial M)$ is the homotopy class of $t_1$ in $\pi_1(S^1)$ in the trivialization of $E|\epsilon$ given by $t_2.$

Let $\Delta$ be a two-simplex of $\text{sk}_2(M).$ Then $\langle\, o(E,t_1)-o(E,t_2), \Delta \,\rangle$ is the difference of the sections $t_1|\partial \Delta$ and $t_2|\partial \Delta$ in a trivialization of $E|\Delta.$ But $t_1$ and $t_2$ agree on $\partial \Delta \setminus \partial M,$ so $$\langle\, o(E,t_1)-o(E,t_2), \Delta \,\rangle = \langle\, o(t_1,t_2), \partial \Delta \cap \partial M \,\rangle.$$ Thus for a relative 2-chain $c \in C_2(M,\partial M)$ we have $$\langle\, o(E,t_1)-o(E,t_2), c \,\rangle = \langle\, o(t_1,t_2), \partial c \cap \partial M\,\rangle,$$ proving that $c_1(\mathfrak{s},t_1)-c_1(\mathfrak{s},t_2) = d(t_1-t_2).$
\end{proof}

\begin{prop} \label{prop:4}
Suppose that $(M,\gamma)$ is a taut sutured manifold, $M$ is open, and the map $\pi_0(A(\gamma))
\to \pi_0(\partial M)$ is surjective. Then $(M,\gamma)$ is balanced.
\end{prop}

\begin{proof}
Since $R_+ = R_+(\gamma)$ and $R_- = R_-(\gamma)$ are both norm minimizing representatives of their homology class in $H_2(M,\gamma)$ and $[R_+] = [R_-],$ we see that $x(R_+) = x(R_-).$ Let $V$ be a component of $R(\gamma).$ Then $V$ is open, so $x(V) = -\chi(V),$ except when $V$ is a disk. Suppose that $V$ is a disk component of say $R_+.$ Then if we push $\partial V$ into $R_-$ we get a curve $C$ in $R_-$ which bounds a disk in $M.$ Since $R_-$ is incompressible, $C$ has to be inessential in $R_-,$ so it bounds a disk in $R_-.$ This argument shows that $R_+$ and $R_-$ have the same number of disk components. Thus $\chi(R_+) = \chi(R_-).$
\end{proof}

\begin{rem} \label{rem:1}
Suppose that $(M,\gamma)$ is a balanced sutured manifold and $H_2(M)
= 0.$ Then $\partial M$ is connected, and so $(M,\gamma)$ is strongly
balanced.

This, together with Proposition \ref{prop:4}, shows that if $(M,\gamma)$ is taut, $M$ is open,
$H_2(M) = 0,$ and $\partial M$ is not a torus which belongs to $\gamma,$ then $(M,\gamma)$ is strongly balanced.
\end{rem}

\begin{defn} \label{defn:7}
Let $S$ be a decomposing surface in a balanced sutured manifold
$(M,\gamma)$ such that the positive unit normal field $\nu_S$ of $S$
is nowhere parallel to $v_0$ along $\partial S.$ This holds for
generic $S.$ We endow $\partial S$ with the boundary orientation.
Let us denote the components of $\partial S$ by $T_1, \dots, T_k.$

Let $w_0$ denote the projection of $v_0$ into $TS,$ this is a
nowhere zero vector field. Moreover, let $f$ be the positive unit
tangent vector field of $\partial S.$ For $1 \le i \le k,$ we define
the \emph{index} $I(T_i)$ to be the number of times $w_0$ rotates
with respect to $f$ as we go around $T_i.$ Then define $$I(S) =
\sum_{i=1}^k I(T_k).$$

Let $p(\nu_S)$ be the projection of $\nu_S$ into $v^{\perp}.$
Observe that $p(\nu_S)|\partial S$ is nowhere zero. For $1 \le i \le
k,$ we define $r(T_i,t)$ to be the rotation of $p(\nu_S)|\partial
T_i$ with respect to the trivialization $t$ as we go around $T_i.$
Moreover, let $$r(S,t) = \sum_{i=1}^k r(T_i,t).$$

We introduce the notation
$$c(S,t) = \chi(S) + I(S) - r(S,t).$$
\end{defn}

The following is \cite[Lemma 3.9]{decomposition}.

\begin{lem} \label{lem:4}
Let $(M,\gamma)$ be a balanced sutured manifold and let $S$ be a
decomposing surface as in Definition \ref{defn:7}.
\begin{enumerate}
\item If $T$ is a component of $\partial S$ such that $T \not\subset
\gamma,$ then $$I(T) = -\frac{|T \cap s(\gamma)|}{2}.$$
\item Suppose that $T_1, \dots, T_a$ are components of $\partial S$
such that $\mathcal{T} = T_1 \cup \dots \cup T_a \subset \gamma$ is
parallel to $s(\gamma)$ and $\nu_S$ points out of $M$ along
$\mathcal{T}.$ Then $I(T_j) = 0$ for $1 \le j \le a;$ moreover,
$$\sum_{j=1}^a r(T_j,t) = \chi(R_+(\gamma)).$$
\end{enumerate}
\end{lem}

\begin{rem} \label{rem:2}
Observe that if the components of the decomposing surface $S$ are $S_1, \dots, S_k,$ and $\nu_S$ is
nowhere parallel to $v_0,$ then
$$c(S,t) = c(S_1,t) + \dots + c(S_k,t).$$
\end{rem}

\begin{defn}
Let $(M, \gamma)$ be a balanced sutured manifold, and let
$(S,\partial S) \subset (M, \partial M)$ be a properly embedded
oriented surface. An element $\mathfrak{s} \in \text{Spin}^c(M,
\gamma)$ is called \emph{outer} with respect to $S$ if there is a
unit vector field $v$ on $M$ whose homology class is $\mathfrak{s}$
and $v_p \neq -(\nu_S)_p$ for every $p \in S.$ Here $\nu_S$ is the
unit normal vector field of $S.$ Let $O_S$ denote the set of outer $\text{Spin}^c$
structures.
\end{defn}

The following is \cite[Lemma 3.10]{decomposition}

\begin{lem} \label{lem:9}
Suppose that $(M,\gamma)$ is a strongly balanced sutured manifold.
Let $t \in T(M,\gamma),$ choose $\mathfrak{s} \in
\text{Spin}^c(M,\gamma)$, and let $S$ be a decomposing surface in
$(M, \gamma)$ as in Definition \ref{defn:7}. Denote the components of $S$
by $S_1, \dots, S_k.$ Then $\mathfrak{s}$ is
outer with respect to $S$ if and only if
\begin{equation}
\left\langle\, c_1(\mathfrak{s},t), [S_i] \,\right\rangle = c(S_i,t) \,\,\text{for every}\,\, 1 \le i \le k.
\end{equation}
In particular, if $\mathfrak{s} \in O_S,$ then
\begin{equation}
\left\langle\, c_1(\mathfrak{s},t), [S] \,\right\rangle = \sum_{i=1}^k c(S_i,t) = c(S,t).
\end{equation}
\end{lem}

\begin{defn}
Suppose that $R$ is a compact, oriented, and open surface. Let $C$
be an oriented simple closed curve in $R.$ If $[C]=0$ in
$H_1(R;\mathbb{Z}),$ then $R \setminus C$ can be written as $R_1 \cup
R_2,$ where $R_1$ is the component of $R \setminus C$ that is
disjoint from $\partial R$ and satisfies $\partial R_1 = C.$ We call
$R_1$ the \emph{interior} and $R_2$ the \emph{exterior} of $C.$

We say that the curve $C$ is \emph{boundary-coherent} if either $[C]
\neq 0$ in $H_1(R;\mathbb{Z}),$ or if $[C]=0$ in $H_1(R;\mathbb{Z})$
and $C$ is oriented as the boundary of its interior.
\end{defn}

\begin{defn} \label{defn:1}
A decomposing surface $S$ in $(M,\gamma)$ is called \emph{nice} if S
is open, $\nu_S$ is nowhere parallel to $v_0,$ and for each component
$V$ of $R(\gamma)$ the set of closed components of $S \cap V$
consists of parallel, coherently oriented, and boundary-coherent simple closed
curves.
\end{defn}

\begin{rem}
Note that every open and groomed decomposing surface becomes nice if we put it into generic position
along the boundary.
\end{rem}

The following theorem is one of the main technical results of the paper \cite{decomposition}, see \cite[Theorem 3.11]{decomposition}.

\begin{thm} \label{thm:4}
Let $(M,\gamma)$ be a strongly balanced sutured manifold;
furthermore, let $(M, \gamma)\rightsquigarrow^S (M', \gamma')$ be a
sutured manifold decomposition along a nice decomposing surface $S.$ Denote the components of $S$ by $S_1, \dots, S_k$ and
choose a trivialization $t \in T(M,\gamma).$ Then
$$SFH(M', \gamma') = \mskip-9mu\bigoplus_{\substack{\mathfrak{s} \in
\text{Spin}^c(M,\gamma) :\\ \langle\,
c_1(\mathfrak{s},t),[S_i]\,\rangle = c(S_i,t)\,\, \forall 1\le i \le k}} \mskip-15mu SFH(M,\gamma,\mathfrak{s}).$$
\end{thm}

\section{Adjunction inequality}

\begin{thm}[Adjunction Inequality] \label{thm:1}
Suppose that the sutured manifold $(M,\gamma)$ is strongly balanced,
and fix a trivialization $t \in T(M,\gamma).$ Let $S \subset M$ be a
nice decomposing surface. If a $\text{Spin}^c$-structure
$\mathfrak{s} \in \text{Spin}^c(M,\gamma)$ satisfies $$\langle\,
c_1(\mathfrak{s},t),[S]\,\rangle < c(S,t),$$ then $SFH(M,\gamma,\mathfrak{s}) =
0.$
\end{thm}

\begin{proof}
Let $S_1, \dots, S_k$ denote the components of $S.$ Then
$$\sum_{i = 1}^k \langle\, c_1(\mathfrak{s},t),[S_i] \,\rangle = \langle\, c_1(\mathfrak{s},t),[S]\,\rangle < c(S,t) = \sum_{i=1}^k c(S_i,t),$$ which implies
that there is an $i$ for which $\langle\, c_1(\mathfrak{s},t),[S_i]\,\rangle < c(S_i,t).$

Now we are going to show that $c(S_i,t) - \langle\, c_1(\mathfrak{s},t),[S_i]\,\rangle$ is
even. From the proof of \cite[Lemma 3.10]{decomposition} it follows that $O_{S_i} \neq \emptyset$ and for $\mathfrak{s}_0 \in O_{S_i}$ we have $\langle\, c_1(\mathfrak{s}_0,t),[S_i]\,\rangle = c(S_i,t).$ Thus $$c(S_i,t) - \langle\, c_1(\mathfrak{s},t),[S_i]\,\rangle = \langle\, c_1(\mathfrak{s}_0,t) - c_1(\mathfrak{s},t), [S_i] \,\rangle = 2\langle\,\mathfrak{s}_0 -\mathfrak{s}, [S] \,\rangle.$$

Add $(c(S_i,t)-\langle\, c_1(\mathfrak{s},t),[S_i]\,\rangle)/2$
compressible 1-handles to $S_i$ to obtain a decomposing surface $S_i'$
in $(M,\gamma).$ Since $S_i$ and $S_i'$ agree in a neighborhood of
$\partial S_i = \partial S_i'$ we have $I(S_i') = I(S_i)$ and $r(S_i',t) =
r(S_i,t).$ Moreover, $$\chi(S_i')=\chi(S_i)-c(S_i,t) + \langle\,
c_1(\mathfrak{s},t), [S_i] \,\rangle.$$ Thus $c(S_i',t) = \langle\,
c_1(\mathfrak{s},t),[S_i]\,\rangle = \langle\,
c_1(\mathfrak{s},t),[S_i']\,\rangle.$ Decomposing $(M,\gamma)$ along $S_i'$ we
get a sutured manifold $(M',\gamma')$ which is not taut since
$R(\gamma')$ is compressible. Thus $SFH(M',\gamma') = 0$ by \cite[Proposition 9.18]{sutured}. Using Theorem \ref{thm:4} $$SFH(M', \gamma') = \mskip-9mu\bigoplus_{\substack{\mathfrak{s'} \in
\text{Spin}^c(M,\gamma) :\\ \langle\,
c_1(\mathfrak{s'},t),[S_i']\,\rangle = c(S_i',t)\,}} \mskip-15mu SFH(M,\gamma,\mathfrak{s'}),$$ so from
$\langle\, c_1(\mathfrak{s},t),[S_i']\,\rangle  = c(S_i',t)$ we get that $SFH(M,\gamma,\mathfrak{s}) =
0.$
\end{proof}

\begin{defn}
Let $(M,\gamma)$ be a balanced sutured manifold. We say that a closed, oriented, one-dimensional submanifold $L \subset \partial M$
is an \emph{$\mathcal{L}$-link} if it is transverse to $v_0^{\perp}|\gamma.$
Given an $\mathcal{L}$-link $L$, a homology class $\alpha \in H_2(M,\partial M)$ such that $\partial \alpha = [L]$ in $H_1(\partial M),$ and a $\text{Spin}^c$-structure $\mathfrak{s} \in \text{Spin}^c(M,\gamma),$ we can define the \emph{rotation number}
$\text{rot}_{\alpha,\mathfrak{s}}(L)$ as follows. Choose a properly embedded, oriented, open surface $S \subset M$ such that
$\partial S = L$ and $[S] = \alpha.$ Furthermore, pick a nowhere zero vector field $v$ on $M$ with $v|\partial M = v_0$ whose homology class is $\mathfrak{s}.$ Then $v^{\perp}|S$ is trivial, let $t_S$ be an arbitrary trivialization. We let
$\text{rot}_{\alpha,\mathfrak{s}}(L)$ be the sum over all components of $L$ of the rotation of $p(\nu_S)$ with respect to $t_S.$
Finally, define $x_{\alpha}(L)$ to be the minimum of $-\chi(S)$ for all surfaces $S$ as above.
\end{defn}

It is straightforward to check that $\text{rot}_{\alpha,\mathfrak{s}}(L)$ is independent of the various choices.
In some sense, the notion of an $\mathcal{L}$-link is analogous to the notion of a Legendrian link in contact topology,
and our rotation number corresponds to the classical rotation number of a Legendrian link. This analogy will be justified
by the following Thurston-Bennequin type inequality.

\begin{cor}
Let $(M,\gamma)$ be a strongly balanced sutured manifold, and suppose that $L \subset \partial M$ is an $\mathcal{L}$-link such that for each component $V$ or $R(\gamma)$ the set of closed components of $L \cap V$ consists of parallel, coherently oriented, nonseparating simple closed curves. If for a $\text{Spin}^c$-structure $\mathfrak{s} \in \text{Spin}^c(M,\gamma)$ we have $SFH(M,\gamma,\mathfrak{s}) \neq 0,$ and $\alpha \in H_2(M,\partial M)$ satisfies $\partial \alpha = [L],$ then
$$x_{\alpha}(L) \ge |\text{rot}_{\alpha,\mathfrak{s}}(L)| - \frac{|L \cap \gamma|}{2}.$$
\end{cor}

\begin{proof}
If $S$ is a properly embedded, oriented, open surface $S \subset M$ such that $\partial S = L,$ then we can perturb $S$ slightly fixing $\partial S$ to get  a nice decomposing surface. So we can apply Theorem \ref{thm:1} to get that
$$\langle\, c_1(S,t), [S] \,\rangle \ge c(S,t) = \chi(S) + I(S) - r(S,t).$$ Since $L$ is transverse to $v_0^{\perp}|\gamma,$ no
component of $L$ lies in $\gamma,$ hence by Lemma \ref{lem:4} we have $I(S) = -|L \cap \gamma|/2.$
Observe that we have three trivializations of $v_0^{\perp}|L,$ namely $t,$ $t_S,$ and $p(\nu_S).$ Furthermore,
$\langle\, c_1(S,t), [S] \,\rangle$ is the rotation of $t$ with respect to $t_S,$ while $r(S,t)$ is the rotation
of $p(\nu_S)$ with respect to $t.$ Hence
$$\langle\, c_1(S,t), [S] \,\rangle + r(S,t) = \text{rot}_{\alpha,\mathfrak{s}}(L),$$ and so
$$-\chi(S) \ge -\text{rot}_{\alpha,\mathfrak{s}}(L) - \frac{|L \cap \gamma|}{2}.$$ Taking the minimum of the left hand
side over all such $S,$ we get
$$x_{\alpha}(L) \ge -\text{rot}_{\alpha,\mathfrak{s}}(L) - \frac{|L \cap \gamma|}{2}.$$ Since
$\text{rot}_{-\alpha,\mathfrak{s}}(-L) = \text{rot}_{\alpha,\mathfrak{s}}(L),$ while $x_{-\alpha}(-L) = x_{\alpha}(L),$ the result follows.
\end{proof}

\begin{rem}
Note that if $\xi$ is a contact structure on $M$ such that $\partial M$ is convex with dividing set $s(\gamma)$ and $L$ is a Legendrian link on $\partial M,$ then the Thurston-Bennequin number of $L$ is precisely $-|L \cap s(\gamma)|/2.$
In the above version of the Thurston-Bennequin inequality, the role of the contact structure $\xi$ is played by the $\text{Spin}^c$-structure $\mathfrak{s},$ and instead of tightness of $\xi$ we have the condition $SFH(M,\gamma,\mathfrak{s}) \neq 0.$ However, if $\xi$ is a tight and $\mathfrak{s}_{\xi}$ is the $\text{Spin}^c$-structure of $\xi,$ we might have
$SFH(M,\gamma,\mathfrak{s}_{\xi}) = 0.$ In the other direction, if $SFH(M,\gamma,\mathfrak{s}) \neq 0,$ it is not clear whether there is a tight contact structure $\xi$ with $\mathfrak{s}_{\xi} = \mathfrak{s}.$
\end{rem}

\begin{thm} \label{thm:3}
Suppose that $(M,\gamma)$ is a strongly balanced sutured manifold and fix a trivialization $t \in T(M,\gamma).$ Let $(M,\gamma) \rightsquigarrow^S (M',\gamma')$ be a surface decomposition along a nice decomposing surface. Then $$SFH(M', \gamma') = \mskip-9mu\bigoplus_{\substack{\mathfrak{s} \in \text{Spin}^c(M,\gamma) :\\ \langle\,
c_1(\mathfrak{s},t),[S]\,\rangle = c(S,t)}} \mskip-15mu SFH(M,\gamma,
\mathfrak{s}).$$
\end{thm}

\begin{proof}
By Theorem \ref{thm:4} we know that if the components of $S$ are $S_1, \dots, S_k$ then $$SFH(M', \gamma') = \mskip-9mu\bigoplus_{\substack{\mathfrak{s} \in
\text{Spin}^c(M,\gamma) :\\ \langle\,
c_1(\mathfrak{s},t),[S_i]\,\rangle = c(S_i,t)\,\, \forall 1\le i \le k}} \mskip-15mu SFH(M,\gamma,\mathfrak{s}).$$ If $\mathfrak{s} \in \text{Spin}^c(M,\gamma)$ satisfies $\langle\, c_1(\mathfrak{s},t),[S_i]\,\rangle = c(S_i,t)$ for every $1\le i \le k,$ then it also satisfies
$\langle\,c_1(\mathfrak{s},t),[S]\,\rangle = c(S,t).$

Now suppose that
\begin{equation} \label{eqn:3}
\langle\,c_1(\mathfrak{s},t),[S]\,\rangle = c(S,t)
\end{equation}
 and $SFH(M,\gamma,\mathfrak{s}) \neq 0.$ If we apply Theorem \ref{thm:1} to $S_i$ we get that $\langle\, c_1(\mathfrak{s},t),[S_i] \,\rangle \ge c(S_i,t).$ Summing these inequalities for $1 \le i \le k$ gives the equality \ref{eqn:3}. Thus $\langle\, c_1(\mathfrak{s},t),[S_i] \,\rangle = c(S_i,t)$ for every $1 \le i \le k.$

So we have shown that $$\mskip-9mu\bigoplus_{\substack{\mathfrak{s} \in
\text{Spin}^c(M,\gamma) :\\ \langle\,
c_1(\mathfrak{s},t),[S_i]\,\rangle = c(S_i,t)\,\, \forall 1\le i \le k}} \mskip-15mu SFH(M,\gamma,\mathfrak{s})= \mskip-9mu\bigoplus_{\substack{\mathfrak{s} \in \text{Spin}^c(M,\gamma) :\\ \langle\,
c_1(\mathfrak{s},t),[S]\,\rangle = c(S,t)}} \mskip-15mu SFH(M,\gamma,
\mathfrak{s}).$$
\end{proof}

%
%
%

\begin{cor} \label{cor:2}
Suppose that $(M,\gamma)$ is a strongly balanced sutured manifold. If $(M,\gamma) \rightsquigarrow^S (M',\gamma')$ is a decomposition along a nice surface such that $[S] = 0$ in $H_2(M,\partial M),$ then either $$SFH(M',\gamma') \cong SFH(M,\gamma),$$ or $SFH(M',\gamma') = 0,$ in which case $(M',\gamma')$ is not taut.
\end{cor}

\begin{proof}
Note that $\langle\, c_1(\mathfrak{s},t),[S]\,\rangle = 0$ since $[S] = 0.$ Then
by Theorem \ref{thm:3} we see that $SFH(M',\gamma') = SFH(M,\gamma)$ if $c(S,t) = 0,$ and $SFH(M',\gamma') = 0$ if $c(S,t) \neq 0.$ By \cite[Theorem 1.4]{decomposition}, the condition $SFH(M',\gamma') = 0$ implies that $(M',\gamma')$ is not taut.
\end{proof}

\begin{cor} \label{cor:5}
Let $(M,\gamma) \rightsquigarrow^S (M',\gamma')$ be a surface decomposition of a strongly balanced sutured manifold. If there is a surface $S'$ disjoint from $S$ such that $S \cup S'$ is nice, $[S \cup S'] = 0$ in $H_2(M,\partial M),$ and $S \cup S'$ gives a taut decomposition, then $SFH(M',\gamma') \cong SFH(M,\gamma).$
\end{cor}

\begin{proof}
Suppose that $S'$ gives a decomposition $(M',\gamma') \rightsquigarrow^{S'} (M'',\gamma'').$ Then
$$SFH(M'',\gamma'') \le SFH(M',\gamma') \le SFH(M,\gamma).$$ If we apply Corollary \ref{cor:2} to $S \cup S',$ then we get that $SFH(M'',\gamma'') \cong SFH(M,\gamma).$ Thus $SFH(M',\gamma') \cong SFH(M,\gamma).$
\end{proof}

\begin{defn} \label{defn:2}
Let $(M,\gamma)$ be a balanced sutured manifold. The \emph{support} of the sutured Floer homology of $(M,\gamma)$ is
$$S(M,\gamma) = \{\, \mathfrak{s} \in \text{Spin}^c(M,\gamma) \,\colon\,
SFH(M,\gamma,\mathfrak{s}) \neq 0 \,\}.$$ Since $SFH(M,\gamma)$ is a
finitely generated Abelian group, $S(M,\gamma)$ is a finite set.
Moreover, if $(M,\gamma)$ is taut, then $S(M,\gamma) \neq \emptyset$ by \cite[Theorem 1.4]{decomposition}.

Let $i \colon H^2(M,\partial M; \mathbb{Z}) \to  H^2(M,\partial M;
\mathbb{R})$ be the map induced by the embedding $\mathbb{Z}
\hookrightarrow \mathbb{R}.$ If $(M,\gamma)$ is strongly balanced
and $t \in T(M,\gamma),$ then we define
$$C(M,\gamma,t) = \{i(c_1(\mathfrak{s},t)) \,\colon\, \mathfrak{s} \in
S(M,\gamma) \} \subset H^2(M,\partial M; \mathbb{R}).$$  Let
$P(M,\gamma,t)$ be the convex hull of $C(M,\gamma,t)$ inside
$H^2(M,\gamma;\mathbb{R}),$ this is a finite polytope. Thus if
$(M,\gamma)$ is taut and $\alpha \in H_2(M,\partial M),$ then
$$c(\alpha,t) = \min \{\,\langle \, c, \alpha\, \rangle \,\colon\, c
\in C(M,\gamma,t) \,\}$$ is a well-defined number.
\end{defn}

\begin{rem} \label{rem:3}
We call $P(M,\gamma,t)$ the \emph{sutured Floer homology polytope} of the sutured manifold $(M,\gamma).$ It follows from Lemma \ref{lem:2} that for $t_1,t_2 \in T(M,\gamma)$ the relationship  $$P(M,\gamma,t_2) = P(M,\gamma,t_1) + i \circ d(t_2 - t_1)$$ holds. So the sutured Floer homology polytope is well defined up to translations in the vector space $H^2(M,\partial M;\mathbb{R})$.
\end{rem}

\begin{cor} \label{cor:1}
Let the sutured manifold $(M,\gamma)$ be taut and  strongly balanced.
Choose a trivialization $t \in T(M,\gamma)$
and let $\alpha \in H_2(M,\partial M)$ be a non-zero element.  Then
\begin{equation} \label{eqn:5}
\max \{\,c(S,t) \,\colon\, S \,\,\text{is a nice
decomposing surface}, [S] = \alpha\,\} \le c(\alpha,t).
\end{equation}
Moreover, for $S$ in the above set $c(S,t) =
c(\alpha,t)$ if and only if $S$ gives a taut decomposition.
If $H_2(M) = 0$ also holds, then inequality \ref{eqn:5} is an equality.
\end{cor}

\begin{proof}
We introduce the notation $$M = \max \{\,c(S,t) \,\colon\, S
\,\,\text{is a nice decomposing surface}, [S] = \alpha\,\}.$$ First
suppose that $S$ is a nice decomposing surface such that $[S] =
\alpha.$ If $c \in C(M,\gamma,t),$ then there is an $\mathfrak{s} \in
S(M,\gamma)$ such that $i(c_1(\mathfrak{s},t)) = c.$ Since
$\mathfrak{s} \in S(M,\gamma),$ we know that
$SFH(M,\gamma,\mathfrak{s}) \neq 0.$ Thus by Theorem \ref{thm:1} we
get that $$c(S,t) \le \langle\, c_1(\mathfrak{s},t),[S] \,\rangle =
\langle\, c, \alpha \, \rangle.$$ Taking the minimum over all $c \in
C(M,\gamma,t),$ we see that $c(S,t) \le c(\alpha, t).$ Since this
holds for every nice $S,$ this implies that $M \le c(\alpha,t).$

Let $(M,\gamma) \rightsquigarrow^S (M',\gamma')$ be a decomposition along
a nice decomposing surface $S$ such that $[S] = \alpha.$ We
saw above that $c(S,t) \le c(\alpha,t).$ Suppose that $c(S,t)
< c(\alpha,t).$ Then for every $\mathfrak{s} \in
\text{Spin}^c(M,\gamma)$ such that $\langle\, c_1(\mathfrak{s},t),
[S] \,\rangle = c(S,t),$ we have $SFH(M,\gamma,\mathfrak{s}) = 0.$
Thus by Theorem \ref{thm:3} we get that
$SFH(M',\gamma') = 0,$ and so \cite[Theorem 1.4]{decomposition}
implies that $(M',\gamma')$ is not taut.

Now assume that $c(S,t) = c(\alpha,t).$ Then by the definition of
$c(\alpha,t)$ there exists an $\mathfrak{s} \in
\text{Spin}^c(M,\gamma)$ such that $\langle\, c_1(\mathfrak{s},t),
[S] \,\rangle = c(\alpha,t)$ and $SFH(M,\gamma,\mathfrak{s}) \neq
0.$ Since $(M,\gamma)$ is taut, it is irreducible, thus
$(M',\gamma')$ is also irreducible. Using Theorem \ref{thm:3} again
we see that $SFH(M',\gamma') \neq 0$ and
so $(M',\gamma')$ is taut by \cite[Proposition 9.18]{sutured}.

Suppose that $H_2(M) = 0.$ By \cite[Lemma 0.7]{Gabai6}, for every $\alpha \neq 0$ in $H_2(M,\partial M)$
there is a nice surface decomposition $(M,\gamma) \rightsquigarrow^S (M',\gamma')$
such that $(M',\gamma')$ is taut and $[S] = \alpha.$ We can assume
that $S$ is open since $H_2(M)=0.$ Using Theorem \ref{thm:3}
$$SFH(M',\gamma') = \bigoplus_{\mathfrak{s} \in \text{Spin}^c(M,\gamma) \colon \langle\,
c_1(\mathfrak{s},t),[S]\,\rangle = c(S,t)}
SFH(M,\gamma,\mathfrak{s}).$$ Since $(M',\gamma')$ is taut, by
\cite[Theorem 1.4]{decomposition} we see that $SFH(M',\gamma') \neq
0.$ Thus there exists an $\mathfrak{s} \in \text{Spin}^c(M,\gamma)$
such that $\langle\, c_1(\mathfrak{s},t), [S] \,\rangle = c(S,t)$
and $SFH(M,\gamma,\mathfrak{s})$ is non-zero, i.e., $\mathfrak{s} \in
S(M,\gamma).$ This implies that $c(\alpha,t) \le M.$ So indeed
we have $c(\alpha,t) = M.$
\end{proof}

\begin{defn} \label{defn:3}
For $\alpha \in H_2(M,\partial M)$ let $$H_{\alpha} = \left \{x \in
H^2(M, \partial M; \mathbb{R}) \colon \langle x, \alpha \rangle \ =
c(\alpha,t) \right \},$$ moreover, $$C_{\alpha}(M,\gamma,t) =
H_{\alpha} \cap C(M,\gamma,t)$$ and similarly
$$P_{\alpha}(M,\gamma,t) = H_{\alpha} \cap P(M,\gamma,t).$$
Finally, we introduce the notation $$SFH_{\alpha}(M,\gamma) = \bigoplus
\left\{SFH(M,\gamma,\mathfrak{s}) \colon i(c_1(\mathfrak{s},t)) \in
C_{\alpha}(M,\gamma,t) \right\}.$$ Note that this is independent of $t$ by
Lemma \ref{lem:2}.
\end{defn}

\begin{prop} \label{prop:1}
Let the sutured manifold $(M,\gamma)$ be taut and strongly balanced.
Fix an element $\alpha \in H_2(M,\partial M).$ Then $P_{\alpha}(M,\gamma,t)$
is the convex hull of $C_{\alpha}(M,\gamma,t)$ and it is a face of the polytope
$P(M,\gamma,t).$ If $S$ is a nice decomposing surface that gives a
taut decomposition $(M,\gamma) \rightsquigarrow^S (M',\gamma')$ and
$[S] = \alpha,$ then
\begin{equation} \label{eqn:1}
SFH(M',\gamma') \cong SFH_{\alpha}(M,\gamma).
\end{equation}
\end{prop}

\begin{proof}
First we prove equation \ref{eqn:1}. If $S$ gives a taut
decomposition, then $c(S,t) = c(\alpha,t)$ by Corollary \ref{cor:1}.
Thus equation \ref{eqn:1} follows from Theorem \ref{thm:3}.

If $\alpha = 0,$ then $c(\alpha,t) = 0.$ So $C_{\alpha}(M,\gamma,t) =
C(M,\gamma,t)$ and $P_{\alpha}(M,\gamma,t) = P(M,\gamma,t),$ hence
Proposition \ref{prop:1} is true for $\alpha = 0.$ Now suppose that
$\alpha \neq 0.$ Then $H_{\alpha}$ is a hyperplane in $H^2(M,
\partial M; \mathbb{R}).$ Using the definition of $c(\alpha,t),$ we
see that $H_{\alpha} \cap C(M,\gamma,t) \neq \emptyset$ and
$\langle\, c, \alpha \,\rangle \ge c(\alpha,t)$ for every $c \in
C(M,\gamma,t).$ Thus $P_{\alpha}(M,\gamma,t)$ is the convex hull of
$C_{\alpha}(M,\gamma,t)$ and is a face of $P(M,\gamma,t).$
\end{proof}

\begin{rem}
It follows from Proposition \ref{prop:1} that if
$i(c_1(\mathfrak{s},t))$ lies in the interior of the polytope
$P(M,\gamma,t),$ then $SFH(M,\gamma,\mathfrak{s})$ dies under any
nice surface decomposition that strictly decreases $SFH(M,\gamma).$
However, we might still be able to obtain information about the interior of the polytope
using decomposing surfaces that are null-homologous in $H_2(M,\partial M).$
\end{rem}

\begin{cor} \label{cor:4}
Let the sutured manifold $(M,\gamma)$ be taut and balanced, and
suppose that $H_2(M)=0.$ Then the following hold.

\begin{enumerate}
\item For every $\alpha \in H_2(M,\partial
M),$ there exists a groomed surface decomposition $(M,\gamma)
\rightsquigarrow^S (M',\gamma')$ such that $(M',\gamma')$ is taut,
$[S] = \alpha,$ and $$SFH(M',\gamma') \cong
SFH_{\alpha}(M,\gamma).$$ If, moreover, $\alpha$ is well groomed, then $S$ can be chosen
to be well groomed.

\item For every  face $F$ of
$P(M,\gamma,t),$ there exists an $\alpha \in H_2(M,\partial M)$ such
that $F = P_{\alpha}(M,\gamma,t).$
\end{enumerate}
\end{cor}

\begin{proof}
First we prove (1). In the case $\alpha =0$ we can choose $S = \emptyset,$ so suppose
that $\alpha \neq 0.$ Then by \cite[Lemma 0.7]{Gabai6} there
exists a groomed surface decomposition $(M,\gamma)
\rightsquigarrow^S (M',\gamma')$ such that $(M',\gamma')$ is taut
and $[S] = \alpha.$ If $\alpha$ is well groomed, then the existence of a well
groomed $S$ follows from Lemma \ref{lem:6}.
Since $H_2(M) = 0,$ we can assume that $S$ has no
closed components, so $S$ is nice. Thus the first part of Corollary
\ref{cor:4} follows from equation \ref{eqn:1}.

Now we prove (2). Let $P = P(M,\gamma,t).$ If $F = P,$ then $\alpha =
0$ works. So suppose that $F$ is a proper face of $P.$
Recall that $P$ is spanned by points lying in the
lattice $L = i(H^2(M,\partial M; \mathbb{Z})).$ Thus there exists an
affine hyperplane of the form $H = H_0 + v_0,$ where $v_0 \in L$ and
$H_0$ is a linear hyperplane spanned by elements of $L,$ and such
that $F = H \cap P.$

Consider the following commutative diagram.
$$\begin{CD}
L @>>> H^2(M,\partial M; \mathbb{R})\\
@VV{u|_L}V @VVuV\\
\text{Hom}(H_2(M,\partial M),\mathbb{Z}) @>>> \text{Hom}(H_2(M,\partial M),\mathbb{R})
\end{CD}$$
Here the horizontal arrows are embeddings. Moreover, $u$ is given by $u(c)(\lambda) = \langle\, c,\lambda\,\rangle$
for $c \in H^2(M,\partial M;\mathbb{R})$ and $\lambda \in H_2(M,\partial M;\mathbb{R}).$
Both $u$ and $u|_L$ are isomorphisms because of the universal coefficient
theorem. Let $b_1,\dots,b_n$ be a basis of the free Abelian group $L.$ Then
$b_1,\dots,b_n$ is also a basis of $H^2(M,\partial M;\mathbb{R}),$ and so defines
a scalar product $\cdot$ on $H^2(M,\partial M;\mathbb{R}).$ Since $u|_L$ is an isomorphism, and because $H_2(M,\partial M)$
is torsion free, there are unique
elements $\beta_1,\dots,\beta_n \in H_2(M,\partial M)$ that satisfy the condition
$\langle\,b_i,\beta_j\,\rangle = \delta_{ij},$ where $\delta_{ij}$ is the Kronecker delta.


Since $H_0$ is spanned by elements of $L,$ there is a vector $a \in
L$ which is perpendicular to $H_0$ and such that $v_0 + a$ and
$P$ lie on the same side of $H.$ In other words,
$a \cdot H_0 = a \cdot (H-v_0) = 0$ and $a \cdot (P-v_0) \ge 0.$
Thus $a \cdot H = a \cdot v_0$ and $a \cdot P \ge a \cdot v_0.$
Let $A_1, \dots, A_n$ be the coordinates of $a$ in the basis $b_1,\dots,b_n,$
these are all integers. Define $\alpha = A_1 \beta_1 + \dots + A_n \beta_n \in H_2(M,\partial M).$
Then for any $c \in H^2(M,\partial M;\mathbb{R}),$ we have $\langle\,c,\alpha \,\rangle = a \cdot c.$
Thus $\langle\, H,\alpha\,\rangle = \langle\,v_0,\alpha\,\rangle$ and
$\langle\, P,\alpha\,\rangle \ge \langle\,v_0,\alpha\,\rangle.$ This implies that
$\langle\, v_0, \alpha \,\rangle = c(\alpha,t),$ and so $H=H_{\alpha}.$ Thus
$F = H \cap P = P_{\alpha}(M,\gamma,t).$
\end{proof}

\begin{prop} \label{prop:3}
Let the sutured manifold $(M,\gamma)$ be taut and balanced, and
suppose that $H_2(M)=0.$ If the polytope $P(M,\gamma,t)$ has $k$
vertices, then $$\text{rk}(SFH(M,\gamma)) \ge k,$$ and there exists a groomed
surface decomposition $(M,\gamma) \rightsquigarrow^S (M',\gamma')$
such that $(M',\gamma')$ is taut and $$\text{rk}(SFH(M',\gamma')) \le
\text{rk}(SFH(M,\gamma))/k.$$
\end{prop}

\begin{proof}
Let $v_1, \dots, v_k$ be the
vertices of $P(M,\gamma,t).$ By Corollary \ref{cor:4}, for every
$1 \le j \le k$ there is a groomed surface decomposition $(M,\gamma)
\rightsquigarrow^{S_j} (M'_j,\gamma'_j)$ such that
$(M'_j,\gamma'_j)$ is taut; furthermore, for $\alpha_j = [S_j]$ we have
$SFH(M'_j,\gamma'_j) = SFH_{\alpha_j}(M,\gamma)$ and
$P_{\alpha_j}(M,\gamma,t) = \{v_j\}.$ By \cite[Theorem
1.4]{decomposition}, we see that $\text{rk}(SFH(M'_j,\gamma'_j)) \ge
1$ for $1 \le j \le k.$ Since the faces $P_{\alpha_j}(M,\gamma,t) = C_{\alpha_j}(M,\gamma,t)$
are pairwise disjoint for $j = 1, \dots,k,$ we get that $$SFH(M,\gamma) \ge
\bigoplus_{j=1}^k SFH_{\alpha_j}(M,\gamma).$$ Thus
$\text{rk}(SFH(M,\gamma)) \ge k,$ and for some $1 \le l \le k$ the inequality
$$\text{rk}(SFH_{\alpha_l}(M,\gamma)) \le \text{rk}(SFH(M,\gamma))/k$$ holds.
So we can choose $S = S_l.$
\end{proof}

\section{How the polytope changes under surface decompositions}

In what follows $b_i(X)$ denotes the $i$-th Betti number of a topological space $X.$

\begin{lem} \label{lem:1}
Suppose that $(M,\gamma)$ is a sutured manifold
and let $(M,\gamma) \rightsquigarrow^S (M',\gamma')$ be a surface
decomposition along the nice decomposing surface $S.$ If $H_2(M) = 0$ or $S=D^2,$ then
$$b_1(M') = b_1(M) + b_1(S) - b_0(S) + b_0(M') - b_0(M).$$ Furthermore, $H_2(M) = 0$
implies that $H_2(M') = 0.$
\end{lem}

\begin{proof}
Let $N(S)$ be a regular neighborhood of $S.$ Since $M' \cap N(S)$
is homotopy equivalent to $S \sqcup S,$
the Mayer-Vietoris sequence of the pair $(M',N(S))$ looks like
\begin{equation} \label{eqn:2}
\dots \to H_j(S) \oplus H_j(S) \to H_j(M') \oplus H_j(S) \to H_j(M) \to \dots
\end{equation}
Note that $H_2(S) = 0$ because $S$ is open. Thus if we write down the sequence \ref{eqn:2}
for $j=2,$ then we see that $H_2(M)=0$ implies $H_2(M')=0.$

Suppose that $H_2(M) =0.$ If we look at the portion of
sequence \ref{eqn:2} starting at $H_2(M)$ and we take the alternating sum of the ranks of
the groups that appear, then we get that
$$2b_1(S) - (b_1(M') + b_1(S)) + b_1(M) -2b_0(S) + (b_0(M')+b_0(S)) - b_0(M) = 0.$$
The result follows. If $S = D^2,$ then $b_2(M) = b_2(M'),$ and we obtain the same conclusion.
\end{proof}

Let us review the relative Maslov grading on Sutured Floer homology,
see \cite[Definiton 8.1]{sutured} and \cite[Definition 8.2]{sutured}.

\begin{defn}
For $\mathfrak{s} \in \text{Spin}^c(M,\gamma),$ the \emph{divisibility} of $\mathfrak{s}$ is
$$\mathfrak{d}(\mathfrak{s}) = \gcd_{\xi \in H_2(M; \mathbb{Z})}
\langle\, c_1(\mathfrak{s}), \xi \,\rangle.$$

Suppose that $\mathfrak{s} \in \text{Spin}^c(M,\gamma),$ and let $(\Sigma,
\boldsymbol{\alpha}, \boldsymbol{\beta})$ be an admissible balanced
diagram for $(M, \gamma).$ Then we define a relative
$\mathbb{Z}_{\mathfrak{d}(\mathfrak{s})}$ grading on $CF(\Sigma,
\boldsymbol{\alpha}, \boldsymbol{\beta}, \mathfrak{s})$ such that
for any $\mathbf{x}, \mathbf{y} \in \mathbb{T}_{\alpha} \cap
\mathbb{T}_{\beta}$ with $\mathfrak{s}(\mathbf{x}) =
\mathfrak{s}(\mathbf{y}) = \mathfrak{s}$ we have
$$\text{gr}(\mathbf{x}, \mathbf{y}) = \mu(\phi) \mod
\mathfrak{d}(\mathfrak{s}),
$$ where $\phi \in \pi_2(\mathbf{x}, \mathbf{y})$ is an arbitrary
homotopy class.
\end{defn}

\begin{defn} \label{defn:4}
Let $(M,\gamma) \rightsquigarrow^S (M',\gamma')$ be a surface decomposition.
If $e \colon M' \hookrightarrow M$ denotes the embedding, then we define
$$F_S = PD \circ e_* \circ (PD')^{-1} \colon H^2(M',\partial M'; \mathbb{R})
\to H^2(M, \partial M; \mathbb{R}),$$ hence the following diagram is commutative.
$$\begin{CD}
H_1(M';\mathbb{R}) @>e_*>> H_1(M;\mathbb{R})\\
@VVPD'V   @VVPDV\\
H^2(M',\partial M'; \mathbb{R}) @>F_S>> H^2(M,\partial M;\mathbb{R})
\end{CD}$$
Here $PD$ and $PD'$ are Poncar\'e duality maps and $e_*$ is the map induced by $e.$
We will use the same symbol $F_S$ to denote the map $H^2(M',\partial M') \to H^2(M,\partial M)$
defined over $\mathbb{Z}$ in a completely analogous way.
\end{defn}

\begin{prop} \label{prop:9}
Let $(M,\gamma) \rightsquigarrow^S (M',\gamma')$ be a nice surface
decomposition of a strongly balanced sutured manifold $(M,\gamma),$
and fix $t \in T(M,\gamma)$ and $t' \in T(M',\gamma').$ Then there are an affine map
$$f_S \colon \text{Spin}^c(M',\gamma') \to \text{Spin}^c(M,\gamma)$$
and an element $c(t,t') \in H^2(M,\partial M;\mathbb{R})$ satisfying the following
three conditions.
\begin{enumerate}
\item
$f_S$ maps onto $O_S,$ and for any $\mathfrak{s} \in O_S$ we have
$$SFH(M,\gamma,\mathfrak{s}) \cong \bigoplus_{\mathfrak{s}' \in
\text{Spin}^c(M',\gamma') \colon f_S(\mathfrak{s}') = \mathfrak{s}}
SFH(M',\gamma',\mathfrak{s}').$$ Furthermore, there is an isomorphism
$$\hat{p} \colon SFH(M',\gamma') \to \bigoplus_{\mathfrak{s}\in O_S} SFH(M,\gamma,\mathfrak{s})$$
such that for every $\mathfrak{s}' \in \text{Spin}^c(M',\gamma')$ and for every $x',y' \in SFH(M',\gamma',\mathfrak{s}')$
we have $\hat{p}(SFH(M',\gamma',\mathfrak{s}')) \subset SFH(M,\gamma,f_S(\mathfrak{s}'))$
and $$\text{gr}(x',y') = \text{gr}(\hat{p}(x'),\hat{p}(y')).$$
\item
If $\mathfrak{s}_1', \mathfrak{s}_2' \in \text{Spin}^c(M',\gamma'),$
then $$F_S(\mathfrak{s}_1' - \mathfrak{s}_2') = f_S(\mathfrak{s}_1') -
f_S(\mathfrak{s}_2') \in H^2(M,\partial M).$$
\item
For every $\mathfrak{s}' \in \text{Spin}^c(M',\gamma')$ we have
$$F_S(i(c_1(\mathfrak{s}',t'))) = i(c_1(f_S(\mathfrak{s}'),t)) + c(t,t').$$
\end{enumerate}
\end{prop}

\begin{proof}
\begin{figure}[b]
\includegraphics{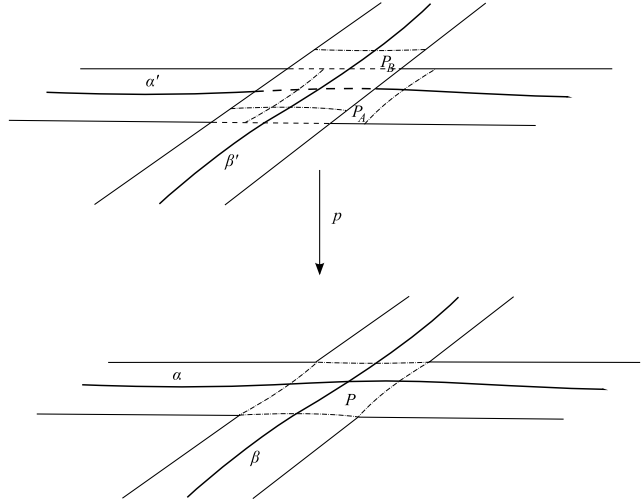}
\caption{The quasi-polygon $P$ in the lower sutured diagram defines a decomposing surface $S.$
The upper sutured diagram corresponds to the manifold obtained by decomposing along $S.$} \label{fig:4}
\end{figure}

We improve the proof of \cite[Theorem 1.3]{decomposition} by also taking into consideration
the $\text{Spin}^c$ and the relative Maslov gradings on $SFH(M',\gamma').$ First we need to recall \cite[Definition 4.3]{decomposition} and \cite[Definition 5.1]{decomposition}.

\begin{defn} \label{defn:16}
A balanced diagram \emph{adapted} to the decomposing surface $S$ in
$(M,\gamma)$ is a quadruple $(\Sigma, \boldsymbol{\alpha},
\boldsymbol{\beta}, P)$ that satisfies the following conditions.
$(\Sigma, \boldsymbol{\alpha}, \boldsymbol{\beta})$ is a balanced
diagram of $(M,\gamma);$ furthermore, $P \subset \Sigma$ is a
quasi-polygon (a closed subsurface of $\Sigma$ whose boundary
is a union of polygons) such that $P \cap \partial
\Sigma$ is exactly the set of vertices of $P.$ We are also given a
decomposition $\partial P = A \cup B,$ where both $A$ and $B$ are
unions of pairwise disjoint edges of $P.$ This decomposition has to
satisfy the property that $\alpha \cap B = \emptyset$ and $\beta
\cap A = \emptyset$ for every $\alpha \in \boldsymbol{\alpha}$ and
$\beta \in \boldsymbol{\beta}.$ Finally, $S$ is given up to
equivalence by smoothing the corners of the surface $(P\times
\{1/2\}) \cup (A \times [1/2,1]) \cup (B \times [0,1/2]) \subset (M,
\gamma).$ The orientation of $S$ is
given by the orientation of $P \subset \Sigma.$ We call a tuple
$(\Sigma, \boldsymbol{\alpha}, \boldsymbol{\beta}, P)$ satisfying
the above conditions a \emph{surface diagram}.
\end{defn}

\begin{defn} \label{defn:18}
Let $(\Sigma, \boldsymbol{\alpha}, \boldsymbol{\beta}, P)$ be a
surface diagram. Then we can uniquely associate to it a tuple
$D(P) = (\Sigma', \boldsymbol{\alpha}',\boldsymbol{\beta}', P_A, P_B, p),$
where $(\Sigma',\boldsymbol{\alpha}', \boldsymbol{\beta}')$ is a balanced diagram,
$p \colon \Sigma' \to \Sigma$ is a smooth map, and $P_A, P_B \subset
\Sigma'$ are two closed subsurfaces (see Figure \ref{fig:4}).

To define $\Sigma',$ take two disjoint copies of $P$ that we call
$P_A$ and $P_B,$ together with diffeomorphisms $p_A \colon P_A \to P$
and $p_B \colon P_B \to P.$ Cut $\Sigma$ along $\partial P$ and
remove $P.$ Then glue $A$ to $P_A$ using $p_A^{-1}$ and $B$ to $P_B$
using $p_B^{-1}$ to obtain $\Sigma'.$ The map $p \colon \Sigma' \to
\Sigma$ agrees with $p_A$ on $P_A$ and with $p_B$ on $P_B,$ and it maps
$\Sigma' \setminus (P_A \cup P_B)$ to $\Sigma \setminus P$ using the
obvious diffeomorphism. Finally, let $\boldsymbol{\alpha}' =
\{\,p^{-1}(\alpha) \setminus P_B \colon \alpha \in
\boldsymbol{\alpha}\,\}$ and $\boldsymbol{\beta}' =
\{\,p^{-1}(\beta) \setminus P_A \colon \beta \in
\boldsymbol{\beta}\,\}.$
\end{defn}

Since $S$ is nice, by \cite[Lemma 4.5]{decomposition} it is isotopic
to a decomposing surface $S'$ such that each component of $\partial
S'$ intersects both $R_+(\gamma)$ and $R_-(\gamma),$ decomposing
$(M,\gamma)$ along $S'$ also gives $(M',\gamma'),$ and $O_S =
O_{S'}.$ Thus we can suppose that each component of $\partial S$
intersects both $R_+(\gamma)$ and $R_-(\gamma).$ Then by
\cite[Proposition 4.8]{decomposition} and \cite[Theorem
6.4]{decomposition} there is a nice and admissible surface diagram
$(\Sigma,\boldsymbol{\alpha},\boldsymbol{\beta},P)$ adapted to $S.$
By \cite[Proposition 5.2]{decomposition}, if $D(P) =
(\Sigma',\boldsymbol{\alpha}',\boldsymbol{\beta}',P_A,P_B,p),$ then
$(\Sigma',\boldsymbol{\alpha}',\boldsymbol{\beta}')$ is an admissible balanced
diagram defining $(M',\gamma').$ Moreover, \cite[Proposition
7.6]{decomposition} says that $p$ gives an isomorphism
$CF(\Sigma',\boldsymbol{\alpha}',\boldsymbol{\beta}') \cong
(O_P,\partial|O_P).$ Here $O_P$ is the subcomplex of
$CF(\Sigma,\boldsymbol{\alpha},\boldsymbol{\beta})$ generated by
$$\left\{\,\mathbf{x} \in \mathbb{T}_{\alpha} \cap \mathbb{T}_{\beta} \colon \mathbf{x} \cap P = \emptyset \,\right\}.$$
However, \cite[Lemma 5.5]{decomposition} implies that $\mathbf{x} \in O_P$ if and only if
$s(\mathbf{x}) \in O_S.$ Thus $p$ induces an isomorphism $$\hat{p} \colon SFH(M',\gamma') \to \bigoplus_{\mathfrak{s}\in O_S} SFH(M,\gamma,\mathfrak{s}).$$

We can now define $f_S.$ If $\mathbb{T}_{\alpha'} \cap \mathbb{T}_{\beta'} \neq \emptyset,$ then fix an element $\mathbf{x}_0' \in \mathbb{T}_{\alpha'} \cap \mathbb{T}_{\beta'}$ and let $\mathbf{x}_0 = p(\mathbf{x}_0')$. Put $\mathfrak{s}_0' = \mathfrak{s}(\mathbf{x}_0')$ and $\mathfrak{s}_0 = \mathfrak{s}(\mathbf{x}_0).$
Otherwise, let $\mathfrak{s}_0' \in \text{Spin}^c(M',\gamma')$ and $\mathfrak{s}_0 \in O_S$ be arbitrary elements. Then for any $\mathfrak{s}' \in \text{Spin}^c(M',\gamma')$ define $$f_S(\mathfrak{s}') = \mathfrak{s}_0 + F_S(\mathfrak{s}'-\mathfrak{s}_0').$$
(2) is immediate from the definition of $f_S.$

Next we show that $f_S$ maps onto $O_S.$ It is sufficient to prove the following claim.

\begin{claim}
Let $\mathfrak{s} \in \text{Spin}^c(M,\gamma).$ Then $\mathfrak{s} \in O_S$ if and only if $\mathfrak{s} - \mathfrak{s}_0 \in \text{im}(F_S).$
\end{claim}

\begin{proof}
Since $\mathbf{x}_0 \in O_P,$ by \cite[Lemma 5.5]{decomposition} we see that $\mathfrak{s}_0 \in O_S.$ Note that $\mathfrak{s} - \mathfrak{s}_0 \in \text{im}(F_S)$ if and only if $PD^{-1}(\mathfrak{s}-\mathfrak{s}_0)$
can be represented by a 1-cycle disjoint from $S.$ Using cut-and-paste techniques, this is equivalent to the statement that
$\langle\,\mathfrak{s}-\mathfrak{s}_0, [S_*]\,\rangle = 0$ for every component $S_*$ of $S.$ Finally, this happens if and only if
$\mathfrak{s}$ and $\mathfrak{s}_0$ can be represented by nowhere zero vector fields that are homotopic (through nowhere zero vector fields) over each component of $S$ rel $\partial S.$ Since $\mathfrak{s}_0 \in O_S,$ this is equivalent to saying that $\mathfrak{s} \in O_S.$
\end{proof}

To prove (1), recall that we have an isomorphism $CF(\Sigma',\boldsymbol{\alpha}',\boldsymbol{\beta}') \cong (O_P,\partial|O_P)$
induced by the projection $p \colon \Sigma' \to \Sigma.$ Moreover, $p(\alpha') = \alpha$ and $p(\beta') = \beta$ for every $\alpha \in \boldsymbol{\alpha}$  and $\beta \in \boldsymbol{\beta}.$ Furthermore, $H_1(M) \cong H_1(\Sigma)/\langle\boldsymbol{\alpha} \cup \boldsymbol{\beta}\rangle$ and $H_1(M') = H_1(\Sigma')/\langle\boldsymbol{\alpha}' \cup \boldsymbol{\beta}'\rangle.$ There is a commutative diagram
$$\begin{CD}
H_1(\Sigma') @>p_*>> H_1(\Sigma)\\
@VV\pi'V   @VV\pi V\\
H_1(M') @>e_*>> H_1(M),
\end{CD}$$
where $\pi$ and $\pi'$ are factor homomorphisms. Now we recall \cite[Definition 4.6]{sutured}.

\begin{defn}
For $\mathbf{x},\mathbf{y} \in \mathbb{T}_{\alpha} \cap \mathbb{T}_{\beta},$ we define
$\epsilon(\mathbf{x}, \mathbf{y}) \in H_1(M)$ as follows. Choose paths $a \colon I \to
\mathbb{T}_{\alpha}$ and $b \colon I \to \mathbb{T}_{\beta}$ with
$\partial a = \partial b = \mathbf{x} - \mathbf{y}.$ Then $a-b$ can be viewed as a one-cycle in $\Sigma$ whose homology class in $M$ is $\epsilon(\mathbf{x}, \mathbf{y}).$ This is independent of the
choices of $a$ and $b.$
\end{defn}

In \cite[Lemma 4.7]{sutured} we showed that $\mathfrak{s}(\mathbf{x}) - \mathfrak{s}(\mathbf{y}) = PD[\epsilon(\mathbf{x},\mathbf{y})]$ for any $\mathbf{x},\mathbf{y} \in \mathbb{T}_{\alpha} \cap \mathbb{T}_{\beta}.$ Now pick elements $\mathbf{x}',\mathbf{y}'
\in \mathbb{T}_{\alpha'} \cap \mathbb{T}_{\beta'},$ and let $\mathbf{x} = p(\mathbf{x}')$ and $\mathbf{y} = p(\mathbf{y}').$
Choose paths $a' \colon I \to \mathbb{T}_{\alpha'}$ and  $b' \colon I \to \mathbb{T}_{\beta'}$ such that $\partial a' = \partial b' =
\mathbf{x}' - \mathbf{y}'.$ Let $a = p(a')$ and $b = p(b').$ Since $p(a'-b') = a-b,$ using the above commutative diagram, we get that
$e_*(\epsilon(\mathbf{x}',\mathbf{y}')) = \epsilon(\mathbf{x},\mathbf{y}).$ Hence
$$PD[e_*(\epsilon(\mathbf{x}',\mathbf{y}'))] = PD[\epsilon(\mathbf{x},\mathbf{y})] = \mathfrak{s}(\mathbf{x}) - \mathfrak{s}(\mathbf{y}).$$ Another application of \cite[Lemma 4.7]{sutured} gives that $$PD[e_*(\epsilon(\mathbf{x}',\mathbf{y}'))] = PD \circ e_* \circ (PD')^{-1}(\mathfrak{s}(\mathbf{x}') - \mathfrak{s}(\mathbf{y}')) = F_S(\mathfrak{s}(\mathbf{x}') - \mathfrak{s}(\mathbf{y}')).$$ So we got that $ F_S(\mathfrak{s}(\mathbf{x}') - \mathfrak{s}(\mathbf{y}')) = \mathfrak{s}(\mathbf{x}) - \mathfrak{s}(\mathbf{y}),$ and by definition $f_S(\mathfrak{s}(\mathbf{x}_0')) = \mathfrak{s}(\mathbf{x}_0).$ Thus $f_S(\mathfrak{s}(\mathbf{x}')) = \mathfrak{s}(\mathbf{x}).$
Hence $$\hat{p}(SFH(M',\gamma',\mathfrak{s}')) \subset SFH(M,\gamma,f_S(\mathfrak{s}'))$$ for every $\mathfrak{s}' \in \text{Spin}^c(M',\gamma').$

If $\mathbf{x}',\mathbf{y}' \in \mathbb{T}_{\alpha'} \cap \mathbb{T}_{\beta'}$ satisfy $\mathfrak{s}(\mathbf{x}') = \mathfrak{s}(\mathbf{y}'),$ then choose a domain $\mathcal{D}'$ connecting $\mathbf{x}'$ and $\mathbf{y}'.$ Then $\mathcal{D} =
p(\mathcal{D}')$ is a domain connecting $\mathbf{x} = p(\mathbf{x}')$ and $\mathbf{y} = p(\mathbf{y}').$ Using Lipshitz's Maslov
index formula (cf. \cite[Proposition 7.3]{decomposition}) we see that $\mu(\mathcal{D}) = \mu(\mathcal{D'})$ since the local diffeomorphism $p$ preserves both the Euler and the point measures. Thus $\text{gr}(x',y') = \text{gr}(\hat{p}(x'),\hat{p}(y'))$
if $x',y' \in SFH(M',\gamma',\mathfrak{s}').$ This concludes the proof of (1).

Finally, we prove (3). Given $t \in T(M,\gamma)$ and $t' \in T(M',\gamma'),$ let
$$c(t,t') = F_S(i(c_1(\mathfrak{s}_0',t'))) - i(c_1(\mathfrak{s}_0,t)).$$
If $\mathfrak{s}' \in \text{Spin}^c(M',\gamma')$ is arbitrary and $\mathfrak{s} =
f_S(\mathfrak{s}'),$ then using (2) we have $$F_S(i(c_1(\mathfrak{s}',t'))) =
F_S(i(c_1(\mathfrak{s}_0',t'))) +
F_S(i(2(\mathfrak{s}'-\mathfrak{s}_0')))=$$ $$ =
i(c_1(\mathfrak{s}_0,t)) + c(t,t') + 2i(\mathfrak{s} - \mathfrak{s}_0) =
i(c_1(\mathfrak{s},t)) + c(t,t').$$
\end{proof}

\begin{rem}
Actually $f_S$ depends only on $S$ and is independent of the choice of a surface diagram representing
$S.$ It can be defined geometrically as follows. Let $v'$ be a vector field representing $\mathfrak{s}' \in \text{Spin}^c(M',\gamma').$
After homotoping $v'$ over $S'_-$ and $S'_+$ (see Definition \ref{defn:8}) one can glue $v'|S'_+$ and $v'|S'_-$ to obtain a nowhere zero vector field $v$ on $M,$ which is unique up to homotopy. Then $v$ represents $f_S(\mathfrak{s}').$ To verify this claim, one has to trace through the identifications in the proof of \cite[Proposition 5.2]{decomposition} to see that $\mathfrak{s}(\mathbf{x}')$ and
$\mathfrak{s}(p(\mathbf{x}'))$ are related by the above "gluing" operation. This is straightforward but tedious, and
we will make no use of it in the rest of the paper.

On the other hand, the isomorphism $\hat{p}$ might possibly depend on the choice of a surface diagram $(\Sigma,\boldsymbol{\alpha},\boldsymbol{\beta},P)$ representing $S.$ To show independence, one either needs a different proof of the decomposition formula, or show invariance of $\hat{p}$ under a sequence of moves relating two different surface diagrams.
\end{rem}

\begin{rem}
If $H^2(M) = 0,$ then for every $t' \in T(M',\gamma')$ there is a $t \in T(M,\gamma)$
such that $c(t,t') = 0$ in Proposition \ref{prop:9}. Indeed, fix an $\mathfrak{s}_0' \in \text{Spin}^c(M',\gamma')$ and let $\mathfrak{s}_0 = f_S(\mathfrak{s}_0').$ The map $H^1(\partial M) \to H^2(M,\partial M)$ is surjective since $H^2(M) = 0.$ Thus by Lemma \ref{lem:2} in the proof of part (3) of Proposition \ref{prop:9}, we can choose a $t
\in T(M,\gamma)$ which satisfies $F_S(i(c_1(\mathfrak{s}_0',t'))) =
i(c_1(\mathfrak{s}_0,t)).$

Note that if we only suppose that $H_2(M) = 0,$ then such a $t$ might not exist. For example, if
$K$ is the knot in $S^1 \times S^2$ which goes around twice monotonically in the $S^1$ direction
and $M = (S^1 \times S^2) \setminus N(K),$ then $M$ fibres over $S^1$ with annulus fibers. Actually $M$ deformation retracts onto a Klein bottle. Thus $H_1(M) = \mathbb{Z}_2 \oplus \mathbb{Z}$ and $H_2(M) = 0.$ So $H^2(M) = \mathbb{Z}_2.$ The map $H^1(\partial M) \to H^2(M,\partial M)$ is Poincar\'e dual to $H_1(\partial M) \to H_1(M).$ This map is $\mathbb{Z} \oplus \mathbb{Z} \to \mathbb{Z}_2 \oplus \mathbb{Z}.$ It takes a pair $(a,b)$ to $(a \mod 2, 2b).$ Here $a$ corresponds to the meridional and $b$ to the longitudinal component. Thus the map $d$ in Lemma \ref{lem:2} is not surjective onto the $\mathbb{Z}$ component of $H^2(M,\partial M).$
\end{rem}

\begin{thm} \label{thm:5}
Let $(M,\gamma) \rightsquigarrow^S (M',\gamma')$ be a taut surface
decomposition, where $(M,\gamma)$ is strongly balanced,
and fix $t \in T(M,\gamma)$ and $t' \in T(M',\gamma').$ Suppose that $S$ is nice and let $\alpha = [S].$ Then the following hold.
\begin{enumerate}
\item[(i)] The map $$F_S - c(t,t') \colon H^2(M',\partial M';\mathbb{R}) \to H^2(M,\partial M;\mathbb{R})$$
projects the polytope $P(M',\gamma',t')$ onto the face $P_{\alpha}(M,\gamma,t)$ of $P(M,\gamma,t).$
\item[(ii)] If, moreover, $S$ is connected and non-separating, then
the image of $F_S$ is the hyperplane $H_{\alpha} + c(t,t').$
\item[(iii)] If, in addition to the assumptions of (2), we have $H_2(M) = 0$ or $S = D^2,$ then $\dim\ker(F_S) = b_1(S).$
\end{enumerate}
\end{thm}

\begin{proof}
Let $\mathfrak{S}= S(M,\gamma)$ and $\mathfrak{S}' = S(M',\gamma').$ First we show that $f_S(\mathfrak{S}') =
\mathfrak{S} \cap O_S.$ Indeed, if $\mathfrak{s}' \in \mathfrak{S}',$ then
by definition $SFH(M',\gamma',\mathfrak{s}') \neq 0.$ So from part (1) of Proposition \ref{prop:9} it
follows that $f_S(\mathfrak{s}') \in O_S$ and $SFH(M,\gamma,f_S(\mathfrak{s}')) \neq 0,$ i.e.,
$f_S(\mathfrak{s}') \in \mathfrak{S} \cap O_S.$ Similarly,
if $\mathfrak{s} \in \mathfrak{S} \cap O_S,$ then by  part (1) of Proposition \ref{prop:9} there exists an $\mathfrak{s}' \in
S(M',\gamma')$ such that $f_S(\mathfrak{s}') = \mathfrak{s}.$

Let $C = C(M,\gamma,t)$ and $C'=C(M',\gamma',t').$ Using part (3) of Proposition \ref{prop:9},
$$F_S(C') - c(t,t') = F_S(i(c_1(\mathfrak{S}',t'))) - c(t,t') = i(c_1(f_S(\mathfrak{S}'),t)) = i(c_1(\mathfrak{S} \cap O_S,t)).$$ By Lemma \ref{lem:9}, for $\mathfrak{s} \in \text{Spin}^c(M,\gamma)$
we have $\mathfrak{s} \in O_S$ if and only if $$\langle\,c_1(\mathfrak{s},t),[S_i]\,\rangle = c(S_i,t)$$
for every component $S_i$ of $S.$ However, we saw in the proof of Theorem \ref{thm:3} that $\mathfrak{s} \in
\mathfrak{S} \cap O_S$ if and only if $\mathfrak{s}
\in \mathfrak{S}$ and $\langle\, c_1(\mathfrak{s},t),[S]\,\rangle = c(S,t).$ Since $S$ gives a taut decomposition,
Corollary \ref{cor:1} gives that $c(S,t) = c(\alpha,t).$ Hence
$$i(c_1(\mathfrak{S} \cap O_S,t)) = C \cap H_{\alpha} = C_{\alpha}(M,\gamma,t),$$ which concludes the proof of (i).

Now we prove (ii). Part (1) of Proposition \ref{prop:9} sates that the image of $f_S$ is
$O_S,$ thus by part (3) the image of $F_S$ is $i(c_1(O_S,t)) + c(t,t').$ Using Lemma \ref{lem:9}, the fact that $S$ is connected,
and $c(S,t) = c(\alpha,t),$ we conclude that $i(c_1(O_S,t)) = H_{\alpha}.$ Since $S$ is non-separating, $H_{\alpha}$ is a hyperplane, i.e., $\dim H_{\alpha} = b_1(M)-1.$

To see (iii), first note that $b_0(M) = b_0(M').$ If $H_2(M) = 0$ or $S = D^2,$ then we can apply Lemma \ref{lem:1} to
conclude that $b_1(M') = b_1(M) + b_1(S) - 1.$ Thus the kernel of the map $F_S$ has dimension $b_1(S).$
\end{proof}

\begin{prop} \label{prop:10}
Suppose that $(M,\gamma) \rightsquigarrow^S (M',\gamma')$ is a taut surface decomposition of a strongly balanced sutured manifold, such that $S$ is a disk and $[S] = \alpha \neq 0.$ Fix $t \in T(M,\gamma)$ and $t' \in T(M',\gamma').$ Then the map $c' \mapsto F_S(c') - c(t,t')$ is an affine isomorphism between the polytope $P(M',\gamma',t')$ and the face $P_{\alpha}(M,\gamma,t)$ of $P(M,\gamma,t).$
\end{prop}

\begin{proof}
We use Theorem \ref{thm:5}. Since $b_1(S) = 0,$ the map $F_S - c(t,t')$ is an affine isomorphism
between $H^2(M',\partial M')$ and $H_{\alpha},$ and maps $P(M',\gamma',t')$ onto $P_{\alpha}(M,\gamma,t).$
\end{proof}

The next result also follows from \cite{Gabai7}.

\begin{cor}
Suppose that the balanced sutured manifold $(M,\gamma)$ is disk decomposable. Then there is a single groomed surface decomposition $(M,\gamma) \rightsquigarrow^S (M',\gamma')$ such that $(M',\gamma')$ is a product. Moreover, if $\gamma$ is connected then $(M,\gamma)$ has a depth at most one taut foliation.
\end{cor}

\begin{proof}
First note that $M$ has to be a handlebody, thus $H_2(M) = 0.$ In this proof we suppress the trivialization $t$ in the notation $P(M,\gamma,t),$ this is justified by Remark \ref{rem:3}. Suppose that $$(M,\gamma) = (M_0,\gamma_0) \rightsquigarrow^{S_1} (M_1,\gamma_1) \rightsquigarrow^{S_2} \dots \rightsquigarrow^{S_n} (M_n,\gamma_n)$$ is a sutured manifold hierarchy such that each $S_i$ is a disk and $(M_n,\gamma_n)$ is a product. Then $SFH(M_k,\gamma_k) \neq 0$ for $0 \le k \le n,$ since $SFH(M_n,\gamma_n) \le SFH(M_k,\gamma_k)$ by Theorem \ref{thm:3}. Furthermore, $(M_k,\gamma_k)$ is irreducible because $(M_n,\gamma_n)$ is. Together with \cite[Proposition 9.18]{sutured} these imply that every $(M_k,\gamma_k)$ is taut.

Let $\alpha_i = [S_i].$ Then $P(M_{i+1},\gamma_{i+1})$ is isomorphic to the face $P_{\alpha_i}(M_i,\gamma_i)$ of $P(M_i,\gamma_i)$ by Proposition \ref{prop:10}. Furthermore, $SFH(M_{i+1},\gamma_{i+1}) \cong SFH_{\alpha_i}(M_i,\gamma_i)$ by Proposition \ref{prop:1}. Since $(M_n,\gamma_n)$ is a product, $P(M_n,\gamma_n)$ is a single point and $SFH(M_n,\gamma_n) \cong \mathbb{Z}.$ So $P(M_n,\gamma_n)$  corresponds to a vertex $v = i(c_1(\mathfrak{s},t))$ of $P(M,\gamma)$ such that $SFH(M,\gamma,\mathfrak{s}) \cong \mathbb{Z}.$

Hence by Corollary \ref{cor:4} there is a groomed and taut surface decomposition $(M,\gamma) \rightsquigarrow^S (M',\gamma')$ such that for $\alpha = [S]$ we have $P_{\alpha}(M,\gamma) = \{v\}$ and $$SFH(M',\gamma') \cong SFH_{\alpha}(M,\gamma) \cong \mathbb{Z}.$$ Since $(M',\gamma')$ is taut, we can use \cite[Theorem 9.7]{decomposition} to conclude that $(M',\gamma')$ is a product. If $\gamma$ is connected, then $S$ is necessarily well groomed, thus by \cite{Gabai} there is a depth at most one taut foliation on $(M,\gamma).$
\end{proof}

\begin{prop} \label{prop:11}
Suppose that $(M,\gamma) \rightsquigarrow^S (M',\gamma')$ is a decomposition of a balanced sutured manifold along a disk $S$ such that $I(S) = -2,$ i.e., $|\partial S \cap s(\gamma)| = 4.$ If we decompose $(M,\gamma)$ along $-S,$ we get $(M',\gamma'').$ Then
\begin{equation} \label{eqn:6}
SFH(M,\gamma) \cong  SFH(M',\gamma') \oplus SFH(M',\gamma'').
\end{equation}
\end{prop}

\begin{proof}
By Theorem \ref{thm:1}, if for some $\mathfrak{s} \in \text{Spin}^c(M,\gamma)$ we have $SFH(M,\gamma,\mathfrak{s}) \neq 0,$ then $\langle\,c_1(\mathfrak{s},t),[S] \,\rangle \ge c(S,t)$ and $\langle\, c_1(\mathfrak{s},t),[-S]\,\rangle \ge c(-S,t).$ Thus $$c(S,t) \le \langle\, c_1(\mathfrak{s},t),[S] \,\rangle \le -c(-S,t).$$ Note that $$-c(-S,t)-c(S,t) = -\chi(S) - \chi(-S) -I(S) - I(-S) = 2.$$ Furthermore, $\langle\, c_1(\mathfrak{s},t),[S] \,\rangle$ is always congruent to $c(S,t)$ modulo $2.$ Indeed, for $\mathfrak{s}_0 \in O_S$ we have $\langle\, c_1(\mathfrak{s}_0,t),[S] \,\rangle = c(S,t)$ and $c_1(\mathfrak{s},t) - c_1(\mathfrak{s}_0,t) = 2(\mathfrak{s} - \mathfrak{s}_0).$

So for any $\mathfrak{s} \in S(M,\gamma)$ either $\langle\, c_1(\mathfrak{s},t),[S] \,\rangle = c(S,t)$ or $\langle\, c_1(\mathfrak{s},t),[-S] \,\rangle = c(-S,t).$ Together with Theorem \ref{thm:3}, this implies equation \ref{eqn:6}.
\end{proof}

\begin{cor} \label{cor:6}
With the assumptions of Proposition \ref{prop:11}, if $(M,\gamma)$ is taut, then at least one of $(M',\gamma')$ and $(M',\gamma'')$ is taut.
\end{cor}

\begin{rem}
Corollary \ref{cor:6} can also be proven using simple cut-and-paste methods. The following, yet unpublished argument was communicated to me by David Gabai. If $x$ denotes the Thurston norm, then $x(R_+(\gamma')) = x(R_+(\gamma)) - 1$ and $x(R_-(\gamma'')) = x(R_-(\gamma))-1.$ If neither $(M',\gamma'),$ nor $(M',\gamma'')$ are taut then we have $x([R_+(\gamma')]) \le x(R_+(\gamma)) -3$ and $x([R_-(\gamma'')]) \le x(R_-(\gamma)) -3.$ Let $T'$ and $T''$ be properly embedded, norm minimizing representatives of $[R_+(\gamma')]$ and $[R_-(\gamma'')],$ respectively. Recall that $S'_+$ and $S'_-$ were introduced in Definition \ref{defn:8}. Let $T$ be the oriented surface obtained from $T' \cup T''$ (viewed as an immersed surface in $M$) by gluing $T' \cap S'_+$ to $T'' \cap S'_-$ and $T' \cap S'_-$ to $T'' \cap S'_+$  (each intersection consists of two arcs), and then doing oriented cut-and-paste along the double curves. We can assume that $T$ has no $S^2$ components since $M$ is irreducible. Then $$x(T) = x([R_+(\gamma')]) + x([R_-(\gamma'')]) + 4 \le x(R_+(\gamma)) + x(R_-(\gamma)) -2 = x(R(\gamma)) -2.$$ Since $[T] = [R(\gamma)]$ in $H_2(M,\gamma),$ we get that $R(\gamma)$ is not norm minimizing in its homology class in $H_2(M,\gamma),$ contradicting the assumption that $(M,\gamma)$ is taut. So at least one of $(M',\gamma')$ and $(M',\gamma'')$ is taut.
\end{rem}

\section{Dimension of the sutured Floer homology polytope}

\begin{thm} \label{thm:2}
Suppose that $H_2(M) = 0$ and the sutured manifold $(M,\gamma)$ is
balanced, taut, reduced and horizontally prime. Let $t \in
T(M,\gamma).$ Then $$\dim P(M,\gamma,t) = \dim H^2(M,\partial M;\mathbb{R}) = b_1(M) = b_1(\partial M)/2.$$
In particular,
$$\text{rk}(SFH(M,\gamma)) \ge b_1(\partial M)/2+1.$$
\end{thm}

\begin{proof}
We improve on and simplify the proof of \cite[Theorem
9.7]{decomposition}. By \cite[Lemma 0.7]{Gabai6} for any non-zero
element $\alpha \in H_2(M,\partial M)$ there is a groomed surface
decomposition $(M,\gamma) \rightsquigarrow^S (M',\gamma')$ such that
$(M',\gamma')$ is taut and $[S] = \alpha.$ We can assume that $S$ is
open since $H_2(M)=0,$ and we can make $S$ nice by putting it into generic position.


We are now going to show that $c(\alpha,t)+ c(-\alpha,t) < 0.$ Since
$S$ gives a taut decomposition, $c(\alpha,t) = c(S,t) = \chi(S) +
I(S) - r(S,t)$ by Corollary \ref{cor:1}. Using \cite[Theorem 2.5]{Scharlemann} and the fact
that $H_2(M) = 0,$ we can find a nice decomposing surface $S'$ which
gives a taut decomposition, $[S'] = -\alpha,$ and $\partial S \cap R(\gamma)
=-\partial S' \cap R(\gamma).$ Then $c(-\alpha,t) = c(S',t) = \chi(S') +
I(S') - r(S',t).$ So we have to show that $$\chi(S)+\chi(S') + I(S)
+ I(S') < r(S,t) + r(S',t).$$ Since $\partial S \cap R(\gamma) = -
\partial S' \cap R(\gamma),$ and by the construction of $S',$ the
one-cycle $\partial S + \partial S' \subset \gamma$ is homologous to
$k s(\gamma)$ in $H_1(\gamma)$ for some non-negative integer $k,$ see Figure 2.

Recall that $r(S,t)$ is defined as the rotation of $p(\nu_S)$ with
respect to the trivialization $t$ as we go around $\partial S,$
where $p$ is orthogonal projection onto $v^{\perp}.$ Note that
$r(S',t)$ can also be computed as the rotation of $-p(\nu_{S'})$
with respect to $t.$ Moreover, $p(\nu_S) = -p(\nu_{S'})$ over $S
\cap R(\gamma) = S' \cap R(\gamma).$

\begin{figure}[t]
\includegraphics{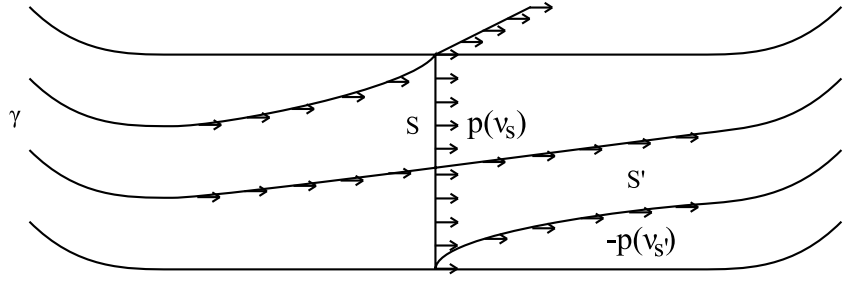}
\caption{The decomposing surfaces $S$ and $S'$} \label{fig:3}
\end{figure}

Let $C$ be a component of
$\partial S \cap \gamma$ or $\partial S' \cap \gamma.$ We can assume
that if $C$ is closed then $\nu_S|C$ (or $-\nu_{S'}|C$) points out of
$M;$ and if $C$ is an arc then it is monotonic between $R_-(\gamma)$
and $R_+(\gamma)$ and $p(\nu_S)|C$ (or $p(\nu_{S'})|C$) is
non-zero and parallel to $s(\gamma).$ Observe that $p(\nu_S) \cup
-p(\nu_{S'})$ is a continuous vector field along $\partial S + \partial
S',$ and it can be homotoped inside $v_0^{\perp} \setminus 0$ such
that it points out of $M$ everywhere. Using the Poncar\'e-Hopf index
formula we get that $$r(S,t) + r(S',t) = k \chi(R_+(\gamma)),$$ see
Lemma \ref{lem:4}.

Suppose that $C$ and $C'$ are components of $\partial S \cap \gamma$
and $\partial S' \cap \gamma,$ respectively, such that $C \approx [0,1] \approx C'$ and
$\partial C = \partial C',$ see Figure \ref{fig:3}.
Let $s_0$ be the component of $s(\gamma)$ containing $C \cap s(\gamma).$
If $C + C'$ is null-homologous in $\gamma,$ then we can achieve by an isotopy of $S'$
that $C = C'.$ If $C + C'$ is homologous to $ms_0$ in $H_1(\gamma)$ for some $m > 0,$
then we can achieve that $|C \cap C'| = m+1.$

Make $S$ and $S'$ transverse by perturbing them in the interior of $M.$
Then take the double curve sum $P$ of $S$ and $S'.$ We are now going to
see how $\chi(S) + \chi(S')$ changes when doing the cut-and-paste. The number of
components of $\partial S \cap \partial S'$ homeomorphic to
$[0,1]$ is at most $|S \cap s(\gamma)|.$ (It is strictly smaller if
there are $[0,1]$ components $C$ and $C'$ of $S \cap \gamma$ and $S' \cap \gamma,$
respectively, such that $C =C'.$) Now let $K$ be a component of $S \cap S'.$
If $\text{int}(K) \subset \text{int}(S) \cap \text{int}(S'),$ then doing
cut-and-paste along $K$ doesn't change the Euler characteristic (we remove
two circles or intervals and glue them back in a different way). On the other
hand, if $K \subset \partial S \cap \partial S'$ is homeomorphic to $[0,1],$
then cut-and-paste along $K$ decreases the Euler characteristic by one
(we glue two surfaces together along two arcs in their boundaries).
According to Lemma \ref{lem:4}, we have $I(S) =
I(S')= -|S \cap s(\gamma)|/2,$ so $$\chi(S) + \chi(S') +
I(S) + I(S') \le \chi(P).$$ Thus it is sufficient to prove that
$\chi(P) < k\chi(R_+(\gamma)).$

Let $r = [R_+(\gamma)] = [R_-(\gamma)] \in H_2(M,\gamma).$ From $H_2(M) = 0,$ and
by looking at the exact sequence of the pair $(M,\gamma),$ we see that the
map $\partial \colon H_2(M,\gamma) \to H_1(\gamma)$ is injective. Thus
$$\partial [P] = \partial (kr) = k [s(\gamma)] \in H_1(\gamma)$$ implies that $[P] = kr$ in
$H_2(M,\gamma).$ Let $x$ denote the Thurston semi-norm on
$H_2(M,\gamma).$ As in \cite[Claim 9.10]{decomposition}, using the
fact that $M$ is irreducible we can suppose that $P$ has no $S^2$
and $T^2$ components.

Suppose that $P$ has a $D^2$ component. Since $H_2(M)=0,$ the
boundary $\partial M$ is connected. Using the fact that $R(\gamma)$
is incompressible and $\partial P \subset \gamma,$ we get that
$\partial M = S^2$ and $\gamma$ is connected. But $M$ is
irreducible, so $M = D^3.$ The sutured manifold $(D^2 \times I,
\partial D^2 \times I)$ obviously satisfies the theorem, so we can
suppose from now on that $P$ has no $D^2$ component. Similarly, we can
assume that $R_+(\gamma)$ has no $D^2$ component.

Using the above assumptions, $x(P) = -\chi(P)$ and $x(R_+(\gamma)) = - \chi(R_+(\gamma)).$
Since $(M,\gamma)$ is taut, $$-\chi(P) = x(P) \ge
x([P]) = x(kr) = kx(r) = -k\chi(R_+(\gamma)).$$ So we only have to
exclude the possibility $x(P) = kx(r).$ In this case, $P$ is norm minimizing in
$kr.$ Thus $P$ cannot have genus $>1$ closed components either because
otherwise we could remove them without changing $[P]$ (as $H_2(M) = 0$) and decrease $x(P).$

Fix a point $z_0 \in R_+(\gamma).$ We define a function $\varphi
\colon M \setminus P \to \mathbb{Z}$ by setting $\varphi(z)$ to be
the algebraic intersection number of $P$ with a path connecting
$z_0$ and $z.$ This is well defined because $[P] = [S]+[S']= \alpha-\alpha = 0$ in
$H_2(M,\partial M),$ and thus any closed curve in $M$ intersects $P$
algebraically zero times. There is a well defined homological
pairing between $H_1(M,R(\gamma))$ and $H_2(M,\gamma).$ Thus if $z
\in R(\gamma),$ then $\varphi(z)$ can be computed by taking the
intersection number of a path connecting $z_0$ and $z$ with the
cycle $kR_-(\gamma).$ So $\varphi|R_+(\gamma) \equiv 0$ and
$\varphi|R_-(\gamma) \equiv k.$ Since $P$ has no closed components,
by considering paths on $\gamma,$ we see that $0 \le \varphi \le k.$

As in \cite[Claim 9.10]{decomposition}, let $J_i =
\text{cl}((\varphi^{-1})(i))$ for $0 \le i \le k$ and let $P_i =
J_{i-1} \cap J_i$ for $1 \le i \le k.$ Then $P$ is the disjoint union of the surfaces
$P_1, \dots, P_k,$ and $\bigcup_{l=0}^{i-1} J_i$ is a homology
between $R_+(\gamma)$ and $P_i$ in $H_2(M,\gamma).$ Thus $[P_i] = r,$
and hence $x(P_i) \ge x(r).$ Since $$\sum_{i=1}^k x(P_i) = x(P) = kx(r),$$
we must have $x(P_i) = x(r)$ for $1 \le i \le k.$ Because $\partial
P$ consists of $k$ parallel copies of $s(\gamma),$ we also see that
$\partial P_i$ is isotopic to $\partial R_+(\gamma)$ for $1 \le i \le k,$ just compute $\varphi|\gamma$
using curves in $\gamma.$ So each $P_i$ is
a horizontal surface. Since $(M,\gamma)$ is horizontally prime, for some $1 \le j \le k$ the surfaces
$P_1, \dots, P_j$ are parallel to $R_+(\gamma)$ and $P_{j+1}, \dots,
P_k$ are parallel to $R_-(\gamma).$

Let $\gamma_j = \gamma \cap J_j.$ Then the sutured manifold
$(J_j,\gamma_j)$ is homeomorphic to $(M,\gamma).$ Thus
$(J_j,\gamma_j)$ is reduced. Observe that the closure of each
component of $S \cap \text{Int}(J_j)$ is either a product disk or product
annulus in $(J_j,\gamma_j)$, which in turn lies in a product
neighborhood $N(\gamma_j)$ of $\gamma_j.$ Indeed, by
Lemma \ref{lem:8} every product disk is inessential in $(M_j,\gamma_j),$
and every product annulus is either ambient isotopic to a component of $\gamma_j$
or bounds a $D^2 \times I$ by Lemma \ref{lem:7}. The rest of $S$, i.e., $S
\setminus J_j$ lies in a product neighborhood of $R(\gamma)$ since
$N_+ = J_1 \cup \dots \cup J_{j-1}$ is a regular neighborhood of
$R_+(\gamma)$ and $N_- = J_{j+1} \cup \dots \cup J_k$ is a regular
neighborhood of $R_-(\gamma).$ Thus $S$ lies in a product
neighborhood $N = N(\gamma) \cup N_+ \cup N_-$ of $\partial M.$ Let
$r \colon N \to \partial M$ be a retraction. Then $r(S)$ represents
a 2-chain in $\partial M$ whose boundary is $\partial S.$ The map
$\partial \colon H_2(M,\partial M) \to H_1(\partial M)$ is injective
because $H_2(M) = 0.$ Moreover, $[S] = \alpha \neq 0,$ thus
$[\partial S] = \partial \alpha \neq 0$ in $H_1(\partial M).$ Consequently,
$\partial S$ cannot be a boundary in $\partial M,$ a contradiction.

So indeed $c(\alpha,t) + c(-\alpha,t) <0.$ By Definition
\ref{defn:2}, this means that the interval $\langle\, \alpha,
P(M,\gamma,t)\,\rangle = [c(\alpha,t), -c(-\alpha,t)]$ is not a
single point. Since this holds for every $\alpha \neq 0$ in $H_2(M,
\partial M),$ the dimensions of $P(M,\gamma,t)$ has to be at least
$b_2(M,\partial M).$ Since $P(M,\gamma,t)$ sits inside
$H^2(M,\partial M, \mathbb{R}),$ the dimension has to be equal to
$b_2(M,\partial M).$ By Poincar\'e duality $b_2(M,\partial M) =
b_1(M),$ and this is equal to $b_1(\partial M)/2$ because
$H_2(M)=0.$

The last statement follows from the fact that a $d$-dimensional
polytope has at least $d+1$ vertices and from Proposition \ref{prop:3}.
\end{proof}

\begin{prop} \label{prop:12}
Suppose that $S$ is a nice decomposing surface in the strongly balanced sutured manifold $(M,\gamma)$ whose components are $S_1, \dots, S_k.$ Let $\alpha = [S]$ and $\alpha_j = [S_j]$ for $1 \le j \le k.$ Suppose that $\alpha \neq 0,$
$$\dim P(M,\gamma,t) = \dim H^2(M,\partial M;\mathbb{R}) = b_1(M),$$ and
$\dim P_{\alpha}(M,\gamma,t) = b_1(M)-1.$  Then there is a non-zero class $\sigma \in H_2(M,\partial M)$ and integers $a_1, \dots, a_k$ such that $\alpha_j = a_j \cdot \sigma$ for $1 \le j \le k.$ Moreover,
$P_{\alpha_j}(M,\gamma,t) = P_{\alpha}(M,\gamma,t)$ for every $1 \le j \le k$ such that $\text{sgn} (a_j) = \text{sgn} (a_1 + \dots +
a_k).$
\end{prop}

\begin{proof}
Note that $\alpha_1 + \dots + \alpha_k = \alpha \neq 0,$ so there is a non-zero $\alpha_j.$ Since $$P_{\alpha}(M,\gamma,t) \subset \bigcap_{j=1}^k H_{\alpha_j},$$ we must have $$\dim \left(\bigcap_{j=1}^k H_{\alpha_j}\right) = b_1(M) - 1.$$ Thus $\alpha_1,\dots,\alpha_k$ are pairwise linearly dependent. The existence of $\sigma$ and $a_1,\dots, a_k$ follows. If $\text{sgn} (a_j) = \text{sgn} (a_1 + \dots + a_k),$ then $\alpha$ and $\alpha_j$ are parallel and point in the same direction, thus $H_{\alpha} = H_{\alpha_j},$ and consequently $P_{\alpha}(M,\gamma,t) = P_{\alpha_j}(M,\gamma,t).$
\end{proof}

\begin{cor}
Let $(M,\gamma)$ be strongly balanced. If $H_2(M) = 0,$ then every face of $P(M,\gamma,t)$ whose dimension is $b_1(M)-1$ is of the form $P_{[R]}(M,\gamma,t)$ for some \emph{connected} groomed surface $R.$
\end{cor}

\begin{proof}
This follows from Proposition \ref{prop:12} and Corollary \ref{cor:4}.
\end{proof}

\section{Depth of a sutured manifold}

\begin{prop} \label{prop:2}
Suppose that $H_2(M) = 0$ and the sutured manifold $(M,\gamma)$ is
balanced, taut, reduced, horizontally prime, and not a product. Then there is always a well-groomed surface decomposition $(M,\gamma) \rightsquigarrow^S (M',\gamma')$ such that
$(M',\gamma')$ is taut and $\text{rk}(SFH(M',\gamma')) \le
\text{rk}(SFH(M,\gamma))/2.$
\end{prop}

\begin{proof}
Note that $\partial M$ is connected because $H_2(M) = 0.$ We show that $\partial M \neq S^2.$ Indeed, otherwise by the irreducibility of $M$ we had $M = D^3,$ and since $(M,\gamma)$ is taut, $\gamma$ had to be a single annulus. But this would contradict that $(M,\gamma)$ is not a product. So $b_1(\partial M) \ge 2.$ By Theorem \ref{thm:2}, $$\dim P(M,\gamma,t) = \dim H^2(M,\partial M;\mathbb{R}) = b_1(\partial M)/2 \ge 1.$$ Lemma \ref{lem:5} implies that there is a well groomed homology class $\alpha \in H_2(M,\partial M).$ Then $-\alpha$ is also well groomed. Thus using Corollary \ref{cor:4} we get well groomed and taut surface decompositions
$(M,\gamma) \rightsquigarrow^{S_1} (M_1,\gamma_1)$ and $(M,\gamma) \rightsquigarrow^{S_2} (M_2,\gamma_2)$ such that $[S_1] = \alpha$ and $[S_2] = -\alpha;$ furthermore, $SFH(M_1,\gamma_1) \cong SFH_{\alpha}(M,\gamma)$ and $SFH(M_2,\gamma_2) \cong SFH_{-\alpha}(M,\gamma).$ Since the dimension of $P(M,\gamma,t)$ is the same as the dimension of the ambient space $H^2(M,\partial M;\mathbb{R}),$ we have $$P_{\alpha}(M,\gamma,t) \cap P_{-\alpha}(M,\gamma,t) = \emptyset.$$ Thus $$SFH_{\alpha}(M,\gamma) \oplus
SFH_{-\alpha}(M,\gamma) \le SFH(M,\gamma),$$ and consequently either $S_1$ or $S_2$ satisfies the requirements of the proposition.
\end{proof}

\begin{rem}
Proposition \ref{prop:2} implies that the number
$\text{rk}(SFH(M,\gamma))$ acts as a complexity of taut balanced
sutured manifolds with $H_2(M) = 0$ in the following sense. If $(M,\gamma)$ is not
a product, then Proposition \ref{prop:15} and Proposition \ref{prop:14} imply that we can perform finitely many horizontal and product annulus decompositions to get an $(M,\gamma)$ which is horizontally prime, reduced, and $H_2(M)$ is still zero. By Proposition \ref{prop:2}, now there is a taut decomposition which strictly decreases $\text{rk}(SFH(M,\gamma)).$

Compare this with the complexity defined in \cite{Gabai} to show the existence of sutured manifold hierarchies. Note that we used the existence of sutured manifold hierarchies to prove that $SFH(M,\gamma) \ge \mathbb{Z}$ if $(M,\gamma)$ is taut, which in turn is implicitly needed in the proof of Proposition \ref{prop:2}.

\end{rem}


\begin{defn} \label{defn:5}
Let $(M,\gamma)$ be a taut sutured manifold. By \cite{Gabai}, $(M,\gamma)$ has a sutured manifold hierarchy $$(M,\gamma) \rightsquigarrow^{S_1} (M_1,\gamma_1) \rightsquigarrow^{S_2} \dots \rightsquigarrow^{S_n} (M_n,\gamma_n),$$ where $(M_n,\gamma_n)$ is a product.
We define the \emph{depth} $d(M,\gamma)$ of $(M,\gamma)$ to be the minimal such $n.$ In particular, $d(M,\gamma) = 0$ if and only if $(M,\gamma)$ is a product.
\end{defn}

\begin{rem}
It is important to note that in the above definition $S_1, \dots, S_n$ can be arbitrary decomposing surfaces, they are not necessarily connected.
\end{rem}

%
%

\begin{figure}[t]
\includegraphics{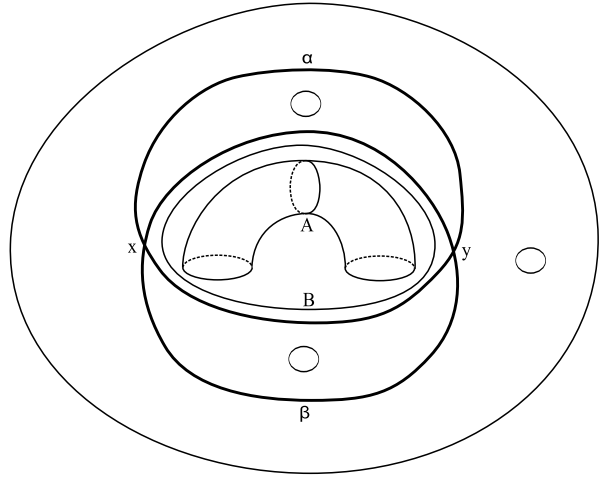}
\caption{A non-product, horizontally prime, taut sutured manifold such that the polytope $P(M,\gamma)$ is a single point.} \label{fig:1}
\end{figure}

\begin{ex} \label{ex:1}
Let $(\Sigma,\boldsymbol{\alpha},\boldsymbol{\beta})$ be the balanced diagram shown in Figure \ref{fig:1}, and let $(M,\gamma)$ be the balanced sutured manifold defined by it. Here $\Sigma$ is a genus one surface with three boundary components, each represented by a little circle. There is one $\alpha$ and one $\beta$ curve; moreover, $\alpha \cap \beta$ consists of two points denoted by $x$ and $y.$ Since there are no periodic domains, $H_2(M) = 0.$

The chain complex $CF(\Sigma,\boldsymbol{\alpha},\boldsymbol{\beta})$ is generated by the points $x$ and $y.$ They lie in the same $\text{Spin}^c$-structure $\mathfrak{s}_0$ because the component $D$ of $\Sigma \setminus (\alpha \cup \beta)$ containing the one-handle gives a homology class of Whitney disks connecting $x$ and $y$ (if we stabilize the diagram there is even a topological Whitney disc in the symmetric product). There are no holomorphic disks connecting $x$ and $y,$ thus $SFH(M,\gamma) \cong \mathbb{Z}^2,$ which lies in $\mathfrak{s}_0.$ Thus $P(M,\gamma,t)$ is a single point. On the other hand, $(M,\gamma)$ is not a product because $SFH(M,\gamma) \ncong \mathbb{Z},$ and it is taut since $SFH(M,\gamma) \neq 0$ and $M$ is irreducible. Moreover, $(M,\gamma)$ is horizontally prime. Indeed, suppose that $$(M,\gamma) \rightsquigarrow^S (M_1,\gamma_1) \sqcup (M_2,\gamma_2)$$ is a horizontal decomposition. Then $$2 = \text{rk}(SFH(M,\gamma)) = \text{rk}(SFH(M_1,\gamma_1)) \cdot \text{rk}(SFH(M_2,\gamma_2)),$$ so $\text{rk}(SFH(M_i,\gamma_i)) = 1$ for $i=1$ or $i=2,$ and this means that $(M_i,\gamma_i)$ is a product. I.e., $S$ is parallel to either $R_+(\gamma)$ or $R_-(\gamma),$ and so $(M,\gamma)$ is horizontally prime. This shows that Theorem \ref{thm:2} fails if $(M,\gamma)$ is not reduced. In fact, there is no nice surface decomposition that would change $SFH(M,\gamma).$ Thus $(M,\gamma)$ cannot be decomposed into a product using a single nice surface decomposition.

Let $A$ denote the core of the handle in $D,$ and let $B$ be a simple closed curve in $D$ parallel to $\partial D.$ Then $A \times I$ is a non-separating and $B \times I$ is a separating product annulus in $(M,\gamma).$ Both of them are nice decomposing surfaces.

If we decompose $(M,\gamma)$ along $A \times I,$ then we get a sutured manifold $(M_A,\gamma_A)$ which is defined by the diagram $(\Sigma_A,\boldsymbol{\alpha},\boldsymbol{\beta}),$ where $\Sigma_A$ is the completion of $\Sigma \setminus A.$ Here $x$ and $y$ lie in different $\text{Spin}^c$-structures. So $SFH(M,\gamma) \cong \mathbb{Z}^2,$ and the two $\mathbb{Z}$-summands lie in different $\text{Spin}^c$-structures. Thus $P(M_A,\gamma_A,t_A)$ consists of two points for any trivialization $t_A.$ By Corollary \ref{cor:4}, there is a well-groomed surface decomposition  $(M_A,\gamma_A) \rightsquigarrow^S (M',\gamma')$ such that $SFH(M',\gamma') \cong \mathbb{Z},$ and thus $(M',\gamma')$ is a product. This shows that $d(M,\gamma) \le 2.$

Decomposing $(M,\gamma)$ along $B \times I,$ we get the disjoint union of a sutured manifold $(M_B,\gamma_B)$ and the product sutured manifold $(D \times I, \partial D \times I).$ Note that $(M_B,\gamma_B)$ can be obtained from $(M_A,\gamma_A)$ by decomposing along a product disk (which corresponds to an arc connecting the feet of the handle in $D$). As above, $SFH(M_B,\gamma_B) \cong \mathbb{Z}^2$ and $P(M_B,\gamma_B,t_B)$ consists of two points. Thus even a separating product annulus can change the sutured Floer homology polytope.

It is not hard to see that $(M_B,\gamma_B)$ is a solid torus with four longitudinal sutures. We can obtain $(M,\gamma)$ from this by attaching $(D \times I, \partial D \times I)$ along $\partial D \times I$ to one of the components of $\gamma_B.$ Of course $D$ is a punctured torus. This again shows that $(M,\gamma)$ is taut. And we can directly see that $(M_B,\gamma_B)$ can be reduced to a product, namely by decomposing along a disk which intersects $s(\gamma_B)$ in four points.
\end{ex}


The following proposition contains \cite[Theorem 9.7]{decomposition}, which claims that $SFH$ detects product sutured manifolds, as the special case $k=0.$ The proof presented here is independent of the proof of \cite[Theorem 9.7]{decomposition}, which refers to an erroneous result in \cite{fibred} that has been corrected in \cite{corrigendum}.

\begin{prop} \label{prop:6}
Suppose that $(M,\gamma)$ is a taut balanced sutured manifold such
that $H_2(M) = 0$ and $\text{rk}(SFH(M,\gamma)) < 2^{k+1}$ for some
integer $k \ge 0.$ Then $$d(M,\gamma) \le 2k.$$
\end{prop}

\begin{proof}
We proceed by induction on $k.$ First suppose that $k=0.$ By Proposition \ref{prop:15}, after a finite sequence of horizontal decompositions we get a taut sutured manifold $(M',\gamma')$ which is horizontally prime. Using Lemma \ref{lem:1}, we see that $H_2(M') = 0$. Furthermore, $SFH(M',\gamma') \cong SFH(M,\gamma)$ by Corollary \ref{cor:2}. Then, using Proposition \ref{prop:14}, we can decompose $(M',\gamma')$ along product annuli to get a reduced, horizontally prime, and taut sutured manifold $(M'',\gamma'').$ Now \cite[Lemma 8.10]{decomposition} and Lemma \ref{lem:1} imply that $SFH(M'',\gamma'') \le SFH(M',\gamma')$ and $H_2(M'') = 0.$ So $\text{rk}(SFH(M'',\gamma'')) \le 1.$ Then, by the second part of Theorem \ref{thm:2}, each component of $\partial M''$ has to be a sphere. But $(M'',\gamma'')$ is taut (in particular irreducible), so it is necessarily a disjoint union of product sutured manifolds of the form $(D^2 \times I, \partial D^2 \times I).$
Consequently, the sutured manifold $(M,\gamma)$ is a product, and hence $d(M,\gamma) = 0.$

Now suppose that $k > 0$ and $(M,\gamma)$ is not a product. First assume that $(M,\gamma)$ is not horizontally prime. Then Proposition \ref{prop:15} gives a  taut decomposition $$(M,\gamma) \rightsquigarrow^H (M',\gamma'),$$ such that each component of $H$ is a horizontal surface and $(M',\gamma')$ is horizontally prime. Let $(M_1',\gamma_1'), \dots, (M_l',\gamma_l')$ denote the components of $(M',\gamma').$ Then $l \ge 2,$ and we can suppose that for each $1 \le i \le l$ the sutured manifold $(M_i',\gamma_i')$ is not a product. Moreover, $H_2(M_i') = 0$ by Lemma \ref{lem:1}. If we apply the $k=0$ case to $(M_i',\gamma_i')$ we get that $\text{rk}(SFH(M_i',\gamma_i')) \ge 2$ for $1 \le i \le l.$ Using Corollary \ref{cor:2} and the K\"unneth formula, we get that $$\text{rk}(SFH(M,\gamma)) = \text{rk}(SFH(M_1',\gamma_1')) \cdot \dots \cdot \text{rk}(SFH(M_l',\gamma_l')).$$ So $\text{rk}(SFH(M_i',\gamma_i')) < 2^k$ for every $1 \le i \le l.$ Hence we can apply the induction hypothesis to each $(M_i',\gamma_i')$ separately to obtain that $d(M_i',\gamma_i') \le 2k-2$ for every $1 \le i \le l.$ But $(M_1',\gamma_1'), \dots, (M_l',\gamma_l')$ are pairwise disjoint, hence $d(M',\gamma') \le 2k-2.$ So $d(M,\gamma) \le 2k-1.$

Consequently, we can assume that $(M,\gamma)$ is horizontally prime. Using Proposition \ref{prop:14}, there is a decomposition $(M,\gamma) \rightsquigarrow^A (M_1,\gamma_1),$ where $A$ is a union of pairwise disjoint product annuli $A_1, \dots, A_r$ and $(M_1,\gamma_1)$ is reduced, horizontally prime, taut, $H_2(M_1) = 0,$ and is not a product. If we apply \cite[Lemma 8.10]{decomposition} to $A_1,\dots,A_r,$ then we get that $SFH(M_1,\gamma_1) \le SFH(M,\gamma).$ So we can use Proposition \ref{prop:3} to get a taut decomposition $(M_1,\gamma_1) \rightsquigarrow^S (M_1',\gamma_1')$ such that $$\text{rk}(SFH(M_1',\gamma_1')) \le
\text{rk}(SFH(M_1,\gamma_1))/2 < 2^k.$$ Thus, using the induction hypothesis on $(M_1',\gamma_1'),$ we see that $d(M_1',\gamma_1') \le 2k-2,$ and hence $d(M,\gamma) \le 2k.$
\end{proof}

\begin{rem} \label{rem:4}
In the above proof, every decomposition can be chosen to be well groomed, except possibly the one along $A,$ which is a disjoint union of product annuli. If every decomposition were well groomed, then we could even claim the existence of a depth at most $2k$ taut foliation on $(M,\gamma).$ Unfortunately, I have overlooked this point in the proof of \cite[Theorem 1.8]{decomposition}, which leaves \cite[Question 9.14]{decomposition} unanswered. If one could make a sutured manifold reduced using a groomed decomposition, that would give a positive answer to \cite[Question 9.14]{decomposition}.
\end{rem}

\begin{cor} \label{cor:7}
Suppose that $K$ is a null-homologous knot in the rational homology 3-sphere $Y,$ and
$$\text{rk}\left(\widehat{HFK}(Y,K,g(K))\right) < 2^{k+1}.$$ Then the sutured manifold $Y(K)$ complementary
to $K$ has depth $d(Y(K)) \le 2k+1.$ In particular, if $k=0,$ then $K$ is fibred.
\end{cor}

\begin{proof}
Let $R$ be a minimal genus Seifert surface for $K.$ By \cite[Theorem 1.5]{decomposition},
$$SFH(Y(R)) \cong \widehat{HFK}(Y,K,g(K)).$$ So Proposition \ref{prop:6} implies that
$d(Y(R)) \le 2k.$ Since we have the sutured manifold decomposition $Y(K) \rightsquigarrow^R Y(R),$
we get $d(Y(K)) \le 2k+1.$ Finally, if $k=0,$ then $Y(R)$ is a product, so $K$ is fibred.
\end{proof}

\section{A semi-norm on the homology of a sutured manifold} \label{sect:7}

In this section, we are going to define a semi-norm on $H_2(M,\partial M;\mathbb{R})$ for a strongly balanced sutured manifold $(M,\gamma).$ Then we will show that it is non-degenerate if $(M,\gamma)$ is taut, reduced, horizontally prime, and $H_2(M) = 0.$ Note that $H_2(M,\partial M)$ is torsion free, and hence can be considered to be a subgroup of $H_2(M,\partial M;\mathbb{R}).$

\begin{defn} \label{defn:6}
Let $(M,\gamma)$ be taut and strongly balanced. For $t \in T(M,\gamma),$ let $p_t \in H^2(M,\partial M;\mathbb{R})$ denote the center of mass of $P(M,\gamma,t).$ Then the polytope
$P(M,\gamma) = P(M,\gamma,t) - p_t$ is independent of $t$ because of Lemma \ref{lem:2}. Of course the center of mass of $P(M,\gamma)$ is $0.$
\end{defn}

\begin{prop} \label{prop:7}
Let $(M,\gamma)$ be taut and strongly balanced. Then for a homology class $\alpha \in H_2(M,\partial M;\mathbb{R})$ the formula  $$y(\alpha) = \max \{\, \langle\,-c,\alpha \,\rangle \colon c \in P(M,\gamma) \,\}$$ defines a semi-norm on $H_2(M,\partial M;\mathbb{R}).$ It is non-degenerate if and only if $$\dim P(M,\gamma) = b_1(M).$$
\end{prop}

\begin{rem}
If $t \in T(M,\gamma),$ then $y(\alpha) = -c(\alpha,t) + \langle\,p_t,\alpha\,\rangle.$
Furthermore, note that for every $k \ge 0$ we have $y(k\alpha) = ky(\alpha),$ but in general $y(\alpha) \neq y(-\alpha)$ can happen.
Indeed, in \cite[Example 8.5]{decateg} we exhibit a family of sutured manifolds whose sutured Floer homology polytopes are all
centrally asymmetric. So in those examples $y$ is not symmetric.
\end{rem}

\begin{proof}
Since $0 \in P(M,\gamma),$ we see that $y(\alpha) \ge 0$ for every $\alpha \in H_2(M,\partial M;\mathbb{R}).$  If $y(\alpha) = 0$ for some $\alpha \neq 0,$ then $P(M,\gamma)$ lies in the hyperplane $$\{c \in H^2(M,\partial M;\mathbb{R}) \colon \langle\,c,\alpha\,\rangle = 0\},$$ thus $$\dim P(M,\gamma) < \dim H^2(M,\partial M;\mathbb{R}) = b_1(M).$$ On the other hand, if $\dim P(M,\gamma) < b_1(M),$ then there is a hyperplane $H$ containing $P(M,\gamma).$ There is also a non-zero homology class $\alpha \in H_2(M,\partial M;\mathbb{R})$ for which $\langle\,H,\alpha\,\rangle = 0,$ i.e., $y(\alpha) = 0$ and $y$ is degenerate.


Suppose that $\alpha,\beta \in H_2(M,\partial M;\mathbb{R}).$ Then $$y(\alpha + \beta) = \max \{\, \langle\,-c,\alpha + \beta \,\rangle \colon c \in P(M,\gamma) \,\} =$$ $$ = \max \{\, \langle\,-c,\alpha \,\rangle + \langle\, -c, \beta \,\rangle \colon c \in P(M,\gamma)\,\} \le$$
$$ \le \max \{\, \langle\,-c,\alpha \,\rangle \colon c \in P(M,\gamma) \,\} + \max \{\, \langle\,-c,\beta \,\rangle \colon c \in P(M,\gamma) \,\} = y(\alpha) + y(\beta).$$
\end{proof}

\begin{rem}
Notice that by construction $-P(M,\gamma)$ is the dual unit norm ball of the semi-norm $y.$
\end{rem}

\begin{prop} \label{prop:8}
Suppose that $(M,\gamma)$ is taut, balanced, reduced, horizontally prime, and $H_2(M) = 0.$ Then
$y$ is a norm on $H_2(M,\partial M;\mathbb{R}).$
\end{prop}

\begin{proof}
Theorem \ref{thm:2} implies that $\dim P(M,\gamma) = b_1(M),$ thus by Proposition \ref{prop:7}
the semi-norm $y$ is non-degenerate.
\end{proof}

\begin{rem}
In \cite{Scharlemann}, another semi-norm is defined on $H_2(M,\partial M;\mathbb{R}),$ which we will denote
by $x^s.$ Given a properly embedded, compact, oriented, and connected surface $S \subset M,$ let
$$x^s(S) = \max \{\,0,-\chi(S) - I(S)\,\},$$ and we extend $x^s$  to disconnected surfaces by taking the
sum over the components. For $\alpha \in H_2(M,\partial M),$ we define $x^s(\alpha)$ as the minimum
of $x^s(S)$ for all properly embedded, compact, oriented surfaces $S$ that represent the homology class $\alpha.$
Finally, it is straightforward to show that $x^s$ extends to $H^2(M,\partial M;\mathbb{R}).$

As opposed to $y,$ the Scharlemann norm $x^s$ is always symmetric. Hence it makes sense to compare $x^s$ with
the symmetrized semi-norm $$z(\alpha) = \frac12(y(\alpha) + y(-\alpha)).$$ In \cite[Theorem 7.7]{decateg},
we show that $z \le x^s.$ Somewhat surprisingly, in general $z \neq x^s$ by \cite[Proposition 7.12]{decateg}.
\end{rem}

\section{Sutured manifolds with $M = S^1 \times D^2$}

In this section, we will compute the sutured Floer homology of every sutured manifold $(M,\gamma)$
with $M = S^1 \times D^2.$ This will illustrate some of the techniques we have just developed.

First note that if such an $(M,\gamma)$ is taut, then $s(\gamma)$ has to be a collection of $n$ parallel
torus knots of type $T_{p,q}.$ Here $p$ denotes the number of times the curve on $\partial M$ goes around in
the longitudinal direction. Furthermore, if $p = 0,$ then $R(\gamma)$ is compressible, hence $(M,\gamma)$ is
not taut. Since $M$ is irreducible, if $(M,\gamma)$ is not taut, then $SFH(M,\gamma) = 0.$ Also note that $n$
is necessarily even. We will denote this sutured manifold by $T(p,q;n).$

\begin{prop}
Let $T(p,q;n)$ be the sutured manifold defined above, and suppose that $n = 2k+2$ for $k \ge 0.$ Then there is an identification $\text{Spin}^c(T(p,q;n)) \cong \mathbb{Z}$ such that the following holds.
\begin{equation} \label{eqn:7}
SFH(T(p,q;n),i) \cong
\begin{cases}
  \mathbb{Z}^{\binom{k}{\lfloor i/p \rfloor}}  & \mbox{if }\,\, 0 \le i < p(k+1), \\
  0 & \mbox{otherwise.}
\end{cases}
 \end{equation}
 Moreover, in each $\text{Spin}^c$-structure any two elements of $SFH$ lie in the same relative Maslov grading.
\end{prop}

\begin{proof}
We saw in Example \ref{ex:1} that $SFH(T(1,0;4)) \cong \mathbb{Z}^2,$ where the two $\mathbb{Z}$ summands lie in $\text{Spin}^c$-structures whose difference is a generator $l$ of $H_1(S^1 \times D^2;\mathbb{Z}).$ Figure \ref{fig:5} shows a sutured diagram for $T(1,0;4).$

\begin{figure}[t]
\includegraphics{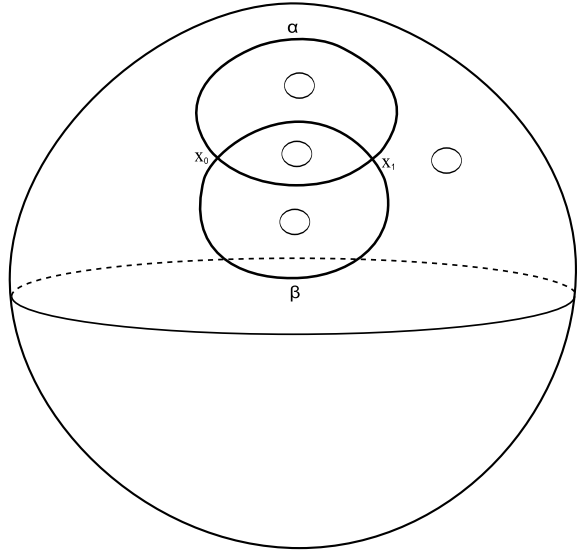}
\caption{A sutured diagram for $T(1,0;4).$} \label{fig:5}
\end{figure}

Let $(M_1,\gamma_1) = T(1,0;n)$ and $(M_2,\gamma_2) = T(p,q;m),$ and suppose that $A_i$ is a component of $\gamma_i$ for $i = 1,2.$ Now glue the annuli $A_1 \subset M_1$ and $A_2 \subset M_2$ such that $A_1 \cap R_-(\gamma_1)$ is identified with $A_2 \cap R_-(\gamma_2).$ Then we obtain the sutured manifold $T(p,q;n+m-2).$ If we decompose $T(p,q;n+m-2)$ along the separating product annulus $A = A_1 = A_2,$ then we get the disjoint union of $T(1,0;n)$ and $T(p,q;m).$ Since $A$ is a nice decomposing surface in $T(p,q;n+m-2),$ an application of Corollary \ref{cor:2} gives that
\begin{equation} \label{eqn:8}
SFH(p,q;n+m-2) \cong SFH(T(1,0;n)) \otimes SFH(T(p,q;m)).
\end{equation}

Using the above formula repeatedly for $T(p,q;m) = T(1,0;4),$ together with the fact that $SFH(T(1,0;4)) \cong \mathbb{Z}^2,$ we get that $$SFH(T(1,0;n)) \cong \bigotimes_{j=1}^k \mathbb{Z}^2,$$ where $n = 2k+2.$ It follows from Proposition \ref{prop:9} that there are generators $x_0^j$ and $x_1^j$ of the $j$-th $\mathbb{Z}^2$ factor in the above expression such that
if $(\eps_1,\dots,\eps_k),(\nu_1,\dots,\nu_k) \in \{0,1\}^k,$ then
$$\mathfrak{s}(x_{\eps_1}^1 \otimes \dots \otimes x_{\eps_k}^k) - \mathfrak{s}(x_{\nu_1}^1 \otimes \dots \otimes x_{\nu_k}^k) = \sum_{j=1}^k (\eps_j -\nu_j) \cdot l.$$ In other words, there is an identification between $\text{Spin}^c(T(1,0;n))$ and $\mathbb{Z}$ such that $\mathfrak{s}(x_{\eps_1}^1 \otimes \dots \otimes x_{\eps_k}^k) = \sum_{j=1}^k \eps_j,$ which proves equation \ref{eqn:7} for $(p,q) = (1,0).$

In light of formula \ref{eqn:8}, we only have to determine $SFH(T(p,q;2)).$ The lens space $L(p,q)$ is obtained
from $(M,\gamma) = T(p,q)$ by gluing an $S^1 \times D^2$ to $M$ such that the meridian $\{1\} \times D^2$ maps to one component, say $\alpha,$ of $s(\gamma).$ Let the knot $K$ be the image of $S^1 \times \{0\}$ in $L(p,q).$ Then the sutured manifold $L(p,q)(K)$ (see \cite[Example 2.4]{sutured}) is exactly $T(p,q;2).$ Hence we only have to find a knot diagram for the knot $K \subset L(p,q).$ But this has already been done in \cite[Proof of Proposition 1.8]{OSz5}. First observe that $K$ is isotopic to a curve on the Heegaard surface $T^2 = S^1 \times S^1$ that intersects $\alpha$ in a single point. Let $\beta$ be a meridian of $M$ that intersects $\alpha$ in exactly $p$ points. Then $(T^2,\alpha,\beta)$ is a Heegaard diagram for $L(p,q).$ As we go around $\alpha,$ label the points
of $\alpha \cap \beta$ with $y_0,\dots,y_{p-1}.$ For $0 \le s \le p-1,$ let $A_s$ be the segment of $\alpha \setminus \{\,y_0,\dots,y_{p-1}\,\}$ connecting $y_s$ and $y_{s+1},$ where $y_p$ is by definition the same as $y_0.$ Choose basepoints $z$ and $w$ on the two sides of $A_{p-1}.$ Then I claim that $(T^2,\alpha,\beta,z,w)$ is a knot diagram defining $K.$ Indeed, if we connect $z$ to $w$ in $T^2 \setminus \alpha$ with an arc, then $w$ to $z$ in $T^2 \setminus \beta$ with a short arc that intersects $\alpha$ in a single point, then we obtain a simple closed curve on $T^2$ that intersects $\alpha$ in a single point, and hence is isotopic to the knot $K.$ Let $\Sigma$ be $T^2$ with two small open disks removed around $z$ and $w.$ Then the previous argument implies that $(\Sigma,\alpha,\beta)$ is a sutured
diagram defining $T(p,q;2).$

It is immediate that $$SFH(T(p,q;2)) \cong \mathbb{Z}^p,$$ which is generated by $y_0, \dots, y_{p-1}.$ Indeed, if we connect $y_s$ and $y_{s+1}$ along $\alpha$ using $A_s$ and then on $\beta$ with an arbitrary curve, then we get a curve on $\Sigma$ whose homology class in $H_1(M)$ is $l$ if $0 \le s <p-1,$ and is $-(p-1)l$ if $s = p$ (this is because components of $\partial \Sigma$ represent $ \pm pl$ in $H_1(M)$). So $\mathfrak{s}(y_{s+1}) - \mathfrak{s}(y_s) = l$ for $0 \le s < p-1.$  This verifies equation \ref{eqn:7} for $n=2$ and $(p,q)$ arbitrary.

To get equation \ref{eqn:7} in general, glue $(M_1,\gamma_1) = T(1,0;n)$ and $(M_2,\gamma_2) = T(p,q;2)$ using formula \ref{eqn:8}, and apply Proposition \ref{prop:9} to see what happens to the $\text{Spin}^c$ grading. If $l_i$ denotes a generator of $H_1(M_i)$ for $i=1,2,$ then $l_1$ is identified with $pl_2$ when we glue $M_1$ and $M_2$ along one of their sutures. This implies equation \ref{eqn:7}.

The above argument actually tells us how to obtain an explicit sutured diagram for $T(p,q;n).$ Let $(\Sigma_1,\alpha_1,\beta_1), \dots, (\Sigma_k,\alpha_k,\beta_k)$ be $k$ identical copies of the sutured diagram
shown in Figure \ref{fig:5}, and introduce the notation $(\Sigma_0,\alpha_0,\beta_0)$ for the diagram defining $T(p,q;2)$ described above. For $0 \le j < k,$ let $c_j$ be a fixed component of $\partial\Sigma_j,$ and for $1 \le h \le k,$ let $d_h$ be a component of $\partial\Sigma_h$ distinct from $c_h.$ Then we obtain $\Sigma$ from $\coprod_{j=0}^k \Sigma_j$ by identifying $c_j$ with $d_{j+1}$ for $0 \le j < k.$ Finally, let $\boldsymbol{\alpha} = \{\,\alpha_0,\alpha_1,\dots\alpha_k\,\}$ and $\boldsymbol{\beta} = \{\,\beta_0,\beta_1,\dots,\beta_k\,\}.$ Then $(\Sigma, \boldsymbol{\alpha},\boldsymbol{\beta})$ is a sutured diagram defining $T(p,q;n).$

Every intersection point in
$\mathbb{T}_{\alpha} \cap \mathbb{T}_{\beta}$ is of the form $y_s \times x_{\eps_1}^1 \times \dots \times x_{\eps_k}^k,$ where $0 \le s \le p-1$ and $(\eps_1,\dots,\eps_k) \in \{0,1\}^k.$ We have $\mathfrak{s}(y_s \times x_{\eps_1}^1 \times \dots \times x_{\eps_k}^k) = \mathfrak{s}(y_{s'} \times x_{\nu_1}^1 \times \dots \times x_{\nu_k}^k)$ if and only if $s = s'$ and $\eps_1 + \dots +\eps_k = \nu_1 + \dots +\nu_k.$ To show that
\begin{equation} \label{eqn:9}
\mu(y_s \times x_{\eps_1}^1 \times \dots \times x_{\eps_k}^k,y_{s} \times x_{\nu_1}^1 \times \dots \times x_{\nu_k}^k) = 0
\end{equation}
for any two intersection points lying in the same $\text{Spin}^c$ structure, it suffices to check the following. If $1 \le t <k$ is fixed and $\nu_j \equiv \eps_j + 1 \mod 2$ for $j=t$ and $j = t+1;$ furthermore,  $\nu_j = \eps_j$ for every other $1 \le j \le k,$ then equation \ref{eqn:9} holds. To see this, look at the region $D$ in $\Sigma \setminus \left(\bigcup\boldsymbol{\alpha} \cup \bigcup \boldsymbol{\beta}\right)$ whose corners are $x_0^t,x_1^t,x_0^{t+1},$ and $x_1^{t+1}.$ This is obtained from the punctured bigons in $\Sigma_t$ and $\Sigma_{t+1}$ containing $c_t$ and $d_{t+1},$ respectively, by gluing $c_t$ to $d_{t+1}.$ Let $\mathcal{D}$ be the domain whose multiplicity in $D$ is one and is zero every where else. Then $\mathcal{D}$ connects $y_s \times x_{\eps_1}^1 \times \dots \times x_{\eps_k}^k$ and $y_s \times x_{\nu_1}^1 \times \dots \times x_{\nu_k}^k;$ moreover, Lipshitz's Maslov index formula \cite[Proposition 7.3]{decomposition} tells us that $\mu(\mathcal{D})= 0.$ This concludes the proof of our claim about the relative Maslov index being zero within a given $\text{Spin}^c$ structure.
\end{proof}

The same way as we proved formula \ref{eqn:8}, we can obtain the following.

\begin{prop}
Suppose that $(M,\gamma)$ is a balanced sutured manifold, and let $\gamma_0$ be a component of $\gamma.$ If $(M,\gamma_1)$ is obtained from $(M,\gamma)$ by adding two sutures parallel to $\gamma_0,$ then
$$SFH(M,\gamma_1) \cong SFH(M,\gamma) \otimes \mathbb{Z}^2.$$
\end{prop}

\begin{rem}
In \cite[Chapter 8]{decateg}, we compute $SFH(M,\gamma,\mathfrak{s})$ for every $\mathfrak{s} \in \text{Spin}^c(M,\gamma)$ when $(M,\gamma)$ is a sutured manifold complementary to a pretzel surface, so $M$ is a genus two handlebody. These examples illustrate well how complicated the support $S(M,\gamma)$ of sutured Floer homology can be in general.
\end{rem}

\bibliographystyle{amsplain}
\bibliography{topology}
\end{document}